\documentclass[11pt,a4]{amsart}
\usepackage{fancyhdr}
\usepackage{amssymb}
\usepackage{amsmath}
\usepackage{array}
\input {xy}

\xyoption{all}
\usepackage[all]{xypic}
\usepackage{latexsym}

\begin{document}

\input{amssym.def}

\newsymbol \circledarrowleft 1309

 %\if@twoside \oddsidemargin 6pt \evensidemargin 6pt \marginparwidth 90pt
 %\else \oddsidemargin 18pt \evensidemargin 18pt \marginparwidth 68pt 
 %\fi
 %\marginparsep 10pt \topmargin -30pt \headheight 12pt \headsep 25pt 
 %\footheight 12pt \footskip 30pt \textheight 680p\textwidth=341ptt \textwidth 455pt 
 %\columnsep 10.5pt \columnseprule 0pt

\newcommand{\ctext}[1]{\makebox(0,0){#1}}
\setlength{\unitlength}{0.1mm}

\newcommand{\wt}{\widetilde}
\newcommand{\Lr}{\Longrightarrow}
\newcommand{\Aut}{\mbox{{\rm Aut}$\,$}}
\newcommand{\ul}{\underline}
\newcommand{\ol}{\overline}
\newcommand{\lr}{\longrightarrow}
\newcommand{\bc}{{\mathbb C}}
\newcommand{\bp}{{\mathbb P}}
\newcommand{\bz}{{\mathbb Z}}

\newcommand{\be}{{\mathbb E}}
\newcommand{\ba}{{\mathbb A}}

\newcommand{\cu}{{\mathcal U}}
\newcommand{\cf}{{\mathcal F}}
\newcommand{\ce}{{\mathcal E}}
\newcommand{\co}{{\mathcal O}}
\newcommand{\cg}{{\mathcal G}}
\newcommand{\cm}{{\mathcal M}}

\newcommand{\mcp}{{\mathcal P}}
\newcommand{\cp}{{\sf P}}

\newcommand{\cq}{{\sf Q}}
\newcommand{\cs}{{\mathcal S}}
\newcommand{\cl}{{\mathcal L}}

\newcommand{\hra}{\hookrightarrow}

\newtheorem{guess}{\sc Theorem}[section]
\newcommand{\bth}{\begin{guess}$\!\!\!${\bf}~~\sl}
\newcommand{\eeth}{\end{guess}}

\newtheorem{propo}[guess]{\sc Proposition}%[section]
\newcommand{\bprop}{\begin{propo}$\!\!\!${\bf }~~\sl}
\newcommand{\eprop}{\end{propo}}
\newtheorem{rema}[guess]{\it Remark}%[section]
\newcommand{\brem}{\begin{rema}$\!\!\!${\bf }~~\rm}
\newcommand{\erem}{\end{rema}}
\newtheorem{coro}[guess]{\sc Corollary}%[section]

\newcommand{\bcor}{\begin{coro}$\!\!\!${\bf }~~\sl}
\newcommand{\ecor}{\end{coro}}
\newtheorem{lema}[guess]{\sc Lemma}%[section]

\newcommand{\blem}{\begin{lema}$\!\!\!${\bf }~~\sl}
\newcommand{\elem}{\end{lema}}

\newtheorem{exam}[guess]{\it Example}%[section]
\newcommand{\beg}{\begin{exam}$\!\!\!${\bf }~~\rm}
\newcommand{\eeg}{\end{exam}}
\newcommand{\er}{\hfill {\Large $\bullet$}\linebreak}

\newtheorem{defe}[guess]{\sc Definition}
\newcommand{\bdefe}{\begin{defe}$\!\!\!${\bf }~~\rm}
\newcommand{\edefe}{\end{defe}}

\newcommand{\spec}{{\rm Spec}\,}

%\pagestyle{myheadings}

%\notitlepage

\title [Donaldson-Uhlenbeck compactification]{Principal bundles
on projective varieties and the Donaldson-Uhlenbeck compactification}
\author{V. Balaji} \address {V.Balaji Chennai Mathematical Institute
92, G.N.Chetty Road, Chennai-600017 INDIA} \email{balaji@cmi.ac.in}

%\date{Preliminary version}

\begin{abstract}
  Let $H$ be a semisimple algebraic group. We prove the semistable
  reduction theorem for $\mu$--semistable principal $H$--bundles over
  a {\it smooth projective variety $X$} defined over the field $\bc$.
  When $X$ is a {\it smooth projective surface} and $H$ is simple, we
  construct the algebro--geometric Donaldson--Uhlenbeck
  compactification of the moduli space of $\mu$--semistable principal
  $H$--bundles with fixed characteristic classes and describe its
  points.  For large characteristic classes we show that the moduli
  space of $\mu$--stable principal $H$--bundles is non--empty.
\end{abstract}
%\vspace{-2mm}

\maketitle

%\vspace{2cm}
%\noindent

%\vspace{2mm}
\section{Introduction}
The purpose of this paper can broadly be termed two-fold. Its first
objective is to prove the {\it semistable reduction theorem} for the
isomorphism classes of $\mu$--semistable principal bundles (in the
sense of Ramanathan--Mumford) with a semisimple structure group $H$
over {\it smooth projective varieties $X$ defined over $\bc$}. In
fact, we prove the semistable reduction theorem for classes of
$\mu$--semistable {\it quasibundles} (Def \ref{quasi}). This
generalises in its entirety, the basic theorem of Langton which proves
that the functor of isomorphism classes of $\mu$--semistable
torsion--free sheaves is {\it proper}. The approach is a
generalisation of the one in \cite{basa} and \cite{bapa}, where this
theorem is proved for curves.

Carrying out the generalisation to the higher dimensional case
involves several new ingredients; for instance, one needs the new
notion of {\it quasibundles} (due to A.Schmitt). This plays the role
of the $\mu$--semistable torsion--free sheaf so as to realise the
boundary points. The final proof is concluded with key inputs from
Bruhat-Tits theory.

Since the proof of the semistable reduction theorem is rather long and
complicated, it is probably appropriate at this point to highlight the
basic differences between our proof and that of Langton in the case of
families of torsion--free sheaves. Except at the very beginning, our
proof follows an entirely different path to that of Langton primarily
because of the fact that it is not even clear if there is a canonical
extension of a family of principal bundles parametrised by a punctured
disc, {\it be it even unstable}, across the puncture. The problem is
no longer sheaf theoretic and one is forced to address the problem of
torsors with structure group which could be {\it non--reductive group
  schemes}. It is to handle this problem that one requires some
aspects of Bruhat-Tits theory. We believe that these new aspects which
come up here should also be of general applicability in similar
situations where compactification questions need to be addressed for
bundles with general structure groups.

More precisely, in his proof Langton first extends the family of
semistable torsion free sheaves to a torsion--free sheaf in the limit
although non-semistable. In other words, the ``structure group'' of
the limiting bundle over a {\it big} open subset still remains $GL_n$.
Then by a sequence of {\it Hecke modifications} he reaches the
semistable limit without changing the isomorphism class of the sheaf
over the generic fibre. Instead, we are forced to extend the family of
semistable {\it rational} $H_K$--bundles to a rational $H_A'$--bundle
with the limiting bundle remaining semistable, but the structure group
becoming {\it non-reductive} in the limit (being the closed fibre of
the group scheme $H_A'$). In other words, one loses the reductivity of
the structure group scheme. Then by using Bruhat-Tits theory, we
relate the group scheme $H_A'$ to the reductive group scheme $H_A$
without changing the isomorphism class of the bundle over the generic
fibre as well as the semistability of the limiting bundle.

Let $A$ be a complete discrete valuation ring and let $K$ be its
quotient field and $k = \bc$ its residue field. Our first main theorem
is the following ($X$ is an arbitrary {\it smooth complex projective
  variety}):
 
\bth\label{}(Theorem \ref{ssred}) Let ${\cp}_K$ be a family of
semistable principal $H$-quasibundles on $X \times \spec K$, or
equivalently, if $H_K$ denotes the group scheme $H \times \spec K$, a
semi-stable $H_K$-quasibundle ${\cp}_K$ on $X_K$. Then there exists a
finite extension $L/K$, with the integral closure $B$ of $A$ in $L$,
such that, ${\cp}_K$, after base change to $\spec B$, extends to a
semistable $H_B$-quasibundle ${\cp}_B$ on $X_B$.  \eeth

We now turn to the second goal of this paper which is to give an
algebro--geometric compactification of the moduli space of
$\mu$--semistable principal bundles over {\it smooth projective
  surfaces}. In fact, we construct a {\it reduced projective scheme}
which can be termed the Donaldson--Uhlenbeck compactification of the
moduli space of $\mu$--stable principal $H$--bundles for a general
simple group $H$ (cf. Corollary \ref{donu}).  In the vector bundle case
such an algebro--geometric construction was given by J.Li (cf.
\cite{li}, \cite{li1}). (see also J.Morgan ( \cite{morgan}).

The theorem of Ramanathan and Subramaniam (\cite{rs}), which is a
generalisation Donaldson-Uhlenbeck-Yau theorem to the case of
principal bundles, gives an identification of antiselfdual (ASD)
Yang-Mills bundles over $X$ with general structure groups with
$\mu$--stable principal bundles. Therefore, our construction of the
Donaldson-Uhlenbeck compactification can be viewed as a natural
compactification of the moduli space of antiselfdual (ASD) Yang-Mills
bundles over $X$ with general structure groups. We remark that for the
case of principal bundles with ASD connections, even a topological
compactification has not been constructed although one can perhaps
extract such a construction from the text \cite{dk}.

More precisely, we have the following theorem:

\bth\label{} (Theorem \ref{uhlenbeck})
\begin{enumerate}
\item Let $H$ be a semisimple algebraic group and $\rho: H \hra SL(V)$
  a faithful rational representation of $H$. There exists a reduced
  projective scheme $M_H(\rho)$ which parametrises equivalence classes
  of $\mu$--semistable $H$--quasibundles with fixed characteristic
  classes, on the {\em smooth projective surface} $X$. This has an
  open subscheme of equivalence classes of $\mu$--semistable principal
  $H$--bundles $M^{0}_H$.

\item Let $H$ be a {\it simple algebraic group}. Then
set--theoretically the closure ${\ol {M^{0}_H}}$ of $\mu$--semistable
principal $H$--bundles $M^{0}_H$ in $M_H(\rho)$ is a subset of the
disjoint union:
\[
{\goth M}_H = \coprod_{{\it l} \geq 0} M_H^{{\mu}-poly}({\sf c}
- {\it l}) \times S^{{\it l}}(X)
\]
where $M_H^{{\mu}-poly}({\sf c} - {\it l})$ is the moduli space of
$\mu$--semistable principal $H$--bundles with characteristic classes
${\sf c} - {\it l}$ (represented as classes of polystable bundles); in
particular, the big stratum is $M_H^{{\mu}-poly}({\sf c}) =
M_H^0$. 

\item When $H$ is {\it simple} the underlying set of points of the
  moduli space ${\ol {M^{0}_H}}$, upto homeomorphism, is independent of
  the representation $\rho$.

\item Let $m = m_{\rho}$ be the Dynkin index of the representation
$\rho$ (Def \ref{dynkin}). The canonical morphism $f_{\rho} :
M_H(\rho) \lr M_{SL(V)}$ maps a copy of $S^{l}(X) \subset M_h(\rho)$
to the symmetric power $S^{m l}(X) \subset M_{SL(V)}$ by sending any
cycle $Z$ to $m \cdot Z$.

\end{enumerate}
\eeth

In the above theorem, the formal construction of the moduli space is
by itself not too difficult. However, the description of its points is
quite involved.  The method of proof is along the lines of the proof
of J.Li (cf. \cite{li}) and the methods in the paper of Le Potier
(\cite{lepot}) (cf. \cite{hl} for a lucid treatment of this approach).
In the description of the points of the moduli space of
$H$--quasibundles and their relationship with the associated moduli
space of $SL(V)$--bundles, the notion of Dynkin index of the
representation $\rho$ makes a natural entry and plays a key role in
defining intrinsically the {\it cycles associated to the points of the
  boundary of the moduli space}. Its significance has already been
noted in the paper by Atiyah, Hitchin and Singer \cite{atiyah} for
bundles on the real four--sphere and in the paper by Kumar, Narasimhan
and Ramanathan \cite{snr} for principal bundles on curves.

Since the construction of the Donaldson--Uhlenbeck compactification
even in the case of vector bundles (\cite{li}) is not entirely by the
methods of GIT it is only natural that for the general case of
arbitrary structure groups which we consider, the use of GIT is only
peripheral. We may recall that the methods of GIT give as a
consequence the projectivity of the quotient space constructed. From
this standpoint, our first theorem (Theorem \ref{ssred}) is absolutely
essential towards proving the compactness. The approach is to separate
the proof of the properness (by proving the semistable reduction), and
the construction of the moduli (by the process of separating points
using sections of a suitable determinant line bundle).

In \cite{donald}, S.Donaldson remarks that it is natural to expect a
generalised theory for his polynomial invariants arising from the
Yang-Mills moduli for bundles with general structure groups. He also
comments that the Uhlenbeck compactness theorem should naturally hold
for the case of general structure groups. One could say that this is
indeed the case in the light of the semistable reduction theorem
mentioned above as well as the description of the points of the moduli
which are added to compactify the ASD moduli space. It is further
remarked in \cite{donald} that if a general theory of these moduli
spaces is given, then one expects results such as the {\it vanishing
  theorems} for these invariants to hold in this general setting. The
existence theorem and the description of points of the moduli in the
present paper aims at securing the foundations of a precise theory
towards this end.

After this paper was completed, the work of Braverman, Finkelberg and
Gaitsgory (\cite{brave}) was brought to our attention by Alexander
Schmitt. They base themselves on a formulation due to V.Ginzburg. The
issues which motivate them in considering Uhlenbeck spaces for
principal bundles with arbitrary structure groups are deep and
far-reaching but quite distinct from ours. Our paper and \cite{brave}
both aim at the construction of these spaces, but the methods used are
altogether different. The coincidence in terminology of {\it
  quasibundles} is also surprising since their notion and ours do not
seem to be related. The approach in \cite{brave} can broadly be termed
{\it ad\`elic} and in the setting of curves this has been used earlier
in \cite{snr}.

We believe that in this paper we have in fact settled affirmatively
some of the basic questions raised in \cite[Page 1]{brave},
particularly those related to {\it moduli}.  It should be very
interesting to establish precise relationships between our paper and
\cite{brave}.

In the final section we prove the following {\it existence} theorem on
which hinges any computation of Donaldson polynomials associated to
these moduli spaces.

\bth\label{} (Theorem \ref{taubes}) Let $H$ be a semisimple algebraic
group over $\bc$. Then the moduli space $M_H(c)^{s}$ of $\mu$--stable
principal $H$--bundles on a smooth projective surface $X$ is
non--empty for {\em large} $c$. \eeth

In the case when $H = SL(2)$ this is highly non-trivial and uses some
deep ideas; this is due to Taubes \cite{taubes} and later due to
Gieseker \cite{gies}. Both methods are deformation theoretic, but the
method used by Taubes is differential geometric ({\it gluing
  techniques}), whereas Gieseker used degeneration techniques in an
algebraic geometric setting to prove the non-emptiness of the moduli
space of $\mu$--stable $SL(2)$--bundles.  Our approach for the general
case of arbitrary semisimple $H$ is to draw on some classical
representation theory, by using what are known as {\it principal
  $SL(2)$'s} in a semisimple group.  We then construct $\mu$--stable
principal $H$--bundles starting from $SL(2)$--bundles for such
principal $SL(2)$'s in $H$. The important point in these existence
results is that the bounds are dependent only on $p_g(X)$ and not on
the polarisation $\Theta$ on $X$. Note that our theorem implies the
non--emptiness of $\mu$--stable $SL(r)$--bundles for all $r$.

The proof of this theorem (where we construct $H$--bundles starting
from $SL(2)$--bundles with $SL(2) \subset H$) and the construction of
the moduli space, where we use a faithful representation $H \subset
SL(V)$, indicate the strong possibility of an {\it algebra} of
Donaldson polynomials coming from the tensor structure on the category
of representations of $H$.

Very recently, G\'omez and Sols (\cite{gomez}) and Schmitt
(\cite{schmitt}, \cite{schmitt1})) have constructed compactifications
of moduli spaces of principal bundles on higher dimensional varieties
using the Gieseker--Maruyama approach for torsion--free sheaves. The
non-emptiness of these spaces (over surfaces) is also therefore a
consequence of our Theorem \ref{taubes}.

G\'omez and Sols follow and generalise the methods of Ramanathan to
higher dimensional varieties. Schmitt gives an alternative approach,
via GIT again, for the moduli construction but in either case this
means that they work with a ``Gieseker--Maruyama'' type definition for
semistability. Schmitt introduces the concept of {\it honest singular
  principal bundles} to recover the boundary points of the moduli
space.  The {\it singular bundles} of Schmitt or equivalently our {\it
  quasibundles} (Def \ref{quasi}) play the key role of giving the
boundary points in our moduli space. It seems possible that the moduli
spaces that we have constructed can be recovered by a generalised
blow-down of the G\'omez--Sols moduli but this needs to be
investigated.

In contrast, the striking feature that emerges here is that the
underlying set of points of our moduli space, ({\it upto
  homeomorphism}), is independent of the choice of a representation of
the structure group, while the moduli spaces of G\'omez-Sols and
Schmitt are invariably dependent on the faithful representation
chosen. We however make no statement on any natural scheme structure
on the moduli space. In fact, this is the case even in the usual
Donaldson-Uhlenbeck compactification for vector bundles.

The brief layout of the paper is as follows: Sections 1 to 4 are
devoted to the proof of the semistable reduction theorem. Section 5 is
devoted to the construction of the compactification and section 6 for
the description of its points. Section 7 contains the proof of the
non--emptiness of the moduli spaces of stable bundles for large
characteristic classes.

\footnotesize{\it Acknowledgments}: I am grateful to C.S.Seshadri and
D.S.Nagaraj for their patience and help in this work.  I thank
M.S.Narasimhan, V.Uma and S.Bandhopadhaya for helpful discussions. I
am extremely grateful to Alexander Schmitt for a meticulous reading of
an earlier version of this paper. His comments have been very helpful
in clarifying many issues in this paper.
%\pagebreak
%\tableofcontents
\normalsize
\begin{center}
{\sc Contents}
\end{center}
{\footnotesize 

\contentsline{section}{\tocsection {}{1}{Introduction}}{}
\contentsline{section}{\tocsection {}{2}{Rational bundles and principal quasibundles}}{}
\contentsline {section}{\tocsection {}{3}{Extension of structure group to the flat closure}}{}
\contentsline {section}{\tocsection {}{4}{Semistable reduction for
quasibundles over projective varieties}}{} 
\contentsline {section}{\tocsection {}{5}{Construction of the moduli
space of bundles over surfaces}}{} 
\contentsline {section}{\tocsection {}{6}{The geometry of the moduli of
$H$--bundles}}{} 
\contentsline {section}{\tocsection {}{7}{Non--emptiness of the moduli
space}}{} 
%\contentsline {section}{\tocsection {}{8}{Appendix}}{} 
\contentsline {section}{\tocsection {}{}{References}}{}
}

\footnotesize
\subsection{Notations and Conventions}
Throughout this paper, unless otherwise stated, we have the
following notations and assumptions:
\small
{\renewcommand{\labelenumi}{{\rm (\alph{enumi})}}
\begin{enumerate}

\item We work over an algebraically closed field $k$ of characteristic zero
and without loss of generality we can take $k$ to be the field of
complex numbers ${\bc}$.

\item {\em $X$ will be a smooth projective variety over $k$ till
\S4. From \S5 onwards it will be a smooth projective surface.}

\item We fix a hyperplane $\Theta$ on $X$ throughout and will use
$\Theta$ for all degree computations.

\item By a {\bf large} or {\bf big} open subset $U \subset X$, we mean
a subset such that $codim_X(X - U) \geq 2$.

\item $H$ is a {\it semisimple} algebraic group, and $G$, unless
otherwise stated will always stand for the general linear group
$GL(V)$ for a finite dimensional vector space $V$.  Their representations
are finite dimensional and rational.

\item $A$ is a discrete valuation ring (which could be assumed to be
complete) with residue field $k$, and quotient field $K$.

\item Let $E$ be a principal $G$-bundle on $X \times T$ where $T$ is
$\spec A$. If $x \in X$ is a closed point then we shall denote by
$E_{x,A}$ or $E_{x,T}$ (resp $E_{x,K}$) the restriction of $E$ to the
subscheme $x \times~\spec A$ or $x \times T$ (resp $x \times~\spec
K$). Similarly, $p \in T$ will denote the closed point of $T$ and the
restriction of $E$ to $X \times p$ will be denoted by $E_p$.

\item In the case of $G = GL(V)$, when we speak of a principal
$G$-bundle we identify it often with the associated vector bundle (and
can therefore talk of the degree of the principal $G$-bundle with
respect to the choice of $\Theta$).

\item We denote by $E_K$ (resp $E_A$) the principal bundle $E$ on $X
\times \spec K$ (resp $X \times \spec A$) when viewed as a principal
$H_K$-bundle (resp $H_A$-bundle). Here $H_K$ and $G_K$ (resp $H_A$ and
$G_A$) are the product group schemes $H \times \spec K$ and $G \times
\spec K$ (resp $H \times \spec A$ and $G \times \spec A$).

\item If $H_A$ is an $A$-group scheme, then by $H_A(A)$ (resp
$H_K(K)$) we mean its $A$ (resp $K$)-valued points. When $H_A = H
\times \spec A$, then we simply write $H(A)$ for its $A$-valued
points. We denote the closed fibre of the group scheme by $H_k$.

\item Let $Y$ be any $G$-variety and let $E$ be a $G$-principal
bundle. For example $Y$ could be a $G$-module. Then we denote by
$E(Y)$ the associated bundle with fibre type $Y$ which is the
following object: $E(Y)$ = $(E \times Y)/G$ for the twisted action of
$G$ on $E \times Y$ given by $g.(e,y)~=~(e.g,g^{-1}.y)$.

\item If we have a group scheme $H_A$ (resp $H_K$) over
$\spec A$ (resp $\spec K$) an $H_A$-module $Y_A$ and a
principal $H_A$-bundle $E_A$. Then we shall denote the
associated bundle with fibre type $Y_A$ by $E_A(Y_A)$.

\item By a family of $H$ bundles on $X$ parametrised by $T$ we mean a
principal $H$-bundle on $X \times T$, which we also denote by
$\{E_t\}_{t \in T}$.
\end{enumerate}

\normalsize
\section{Rational bundles and principal quasibundles}
Let $X$ be a smooth projective variety over $\bc$.  Let $\rho: H \hra
SL(V)$. Let $\ce$ be a torsion--free sheaf on $X$ and let $U({\ce})$ be
the largest open subset where $\ce$ is locally free. Let $S_{\ce}$ be
the affine $X$--scheme given by $\underline {Spec}~Sym^*(\ce \otimes
{\co}(V))$. There is a canonical action of $H$ on $S_{\ce}$ and we
consider the categorical quotient $\underline {Spec}~(Sym^*(\ce
\otimes {\co}(V))^{H}) = (S_{\ce})//H$. Since $H$ is assumed to be
semisimple, the action has a non-empty collection of semistable
points.

Suppose that we are given a morphism ({\it which we term a reduction
section for the obvious reasons}) $\sigma : X \lr (S_{\ce})//H $,
which on the open subset $U(\ce)$ gives a {\it genuine reduction of
structure group} to $H$. We obtain an $X$--scheme $\cp \lr X$, by
pulling back the quotient map $q : S_{\ce} \lr (S_{\ce})//H $. This
notion is due to Alexander Schmitt, and the pair $(\ce, \sigma)$ is
termed by him an {\it honest singular} principal bundle. This is a
natural generalisation of the classical notion of a {\it frame bundle}
associated to a vector bundle.

\bdefe\label{quasi}{\it (Rational principal bundles and principal
quasibundles)}
\begin{enumerate}

\item Following (\cite{rr}), by a rational principal $H$--bundle, we
mean a principal $H$--bundle on a large open $U \subset X$.

\item (A.Schmitt) By an $H$--quasibundle, we mean a scheme $\cp
\lr X$ as above. Let $U(\cp)$ denote the largest open subset of $X$
where $\cp$ is a genuine principal $H$--bundle. In particular,
${\cp}|_{U(\cp)}$ is a rational principal $H$--bundle.
\end{enumerate}
\edefe

\brem\label{quasidiag} We remark that whenever we speak of a
quasibundle we have the following collection of objects:
\begin{enumerate}
\item A faithful representation $H \hra SL(V)$.

\item A torsion--free sheaf $\ce$ with generic fibre of type $V^*$.

\item A diagram:

\[
\begin{array}{ccc}
{\cp} &
\stackrel{}{\longrightarrow} &
S_{\ce} \\
\Big\downarrow & &
\Big\downarrow \\
X   & \stackrel{\sigma}{\longrightarrow} &
S_{\ce}//H
\end{array}
\]

\end{enumerate}
\erem

\brem \label{bigopen} Since the notion of quasibundles will play a key
role in what proceeds, we will briefly recall its salient features,
especially those which will be frequently used in this paper. 

Let $T$ be an arbitrary normal variety and $\ce$ a torsion--free sheaf
on $T$. We can identify the affine $T$--scheme $S_{\ce}$ with the
space ${\underline {Hom}}~~({\ce}, V^{*} \otimes {\co}_X)$ and
similarly the affine $T$--scheme $S_{\ce}//H$ with the space
${\underline {Hom}}~~({\ce}, V^{*} \otimes {\co}_X)//H$ (cf.
\cite[3.7,3.8]{schmitt}).

Let $\cu$ be the maximal open subset of $T$ such that
${\ce}|_{\cu}$ is locally free with general fibre $V$ and trivial
determinant and let $\rho: H \hra GL(V)$. Then one knows that $\cu$ is
a {\it big} open subset of $T$. A reduction of structure group of the
principal $GL(V)$--bundle underlying the vector bundle ${\ce}|_{\cu}$
can be viewed as a section:
\[
{\sigma}_{\cu} : {\cu} \lr {\underline {Isom}}~~({\ce}|_{\cu}, V
\otimes {\co}_X)/H
\]
Now observe that ${\underline {Isom}}~~({\ce}|_{\cu}, V \otimes
{\co}_X)/H \subset S_{\ce}//H$ and further, since $S_{\ce}//H$ can be
{\it embedded as closed subscheme of a vector bundle over $T$}, it
follows by Hartogs theorem that the reduction ${\sigma}_{\cu}$
extends uniquely to a section:
\[
\sigma: X \lr S_{\ce}//H.
\]
\erem

Following \cite{r1} we have the following definitions.

\bdefe\label{ramaa} (A. Ramanathan) A rational principal $H$--bundle
$E$ is said to be {\it ${\mu}$--semistable} (resp. {\it
${\mu}$--stable}) if $\forall$ parabolic subgroup $P$ of $H$,
$\forall$ reduction $\sigma_P : X \lr {\cp}(H/P)$ and $\forall$
dominant character $\chi$ of $P$, the bundle $\sigma_P^* (L_\chi))$
has degree $\leq 0$ (resp. $< 0$).(cf.\cite {r1}). We note that in
this convention, a dominant character $\chi$ of $P$ induces a negative
ample line bundle on $G/P$. Note further that this definition makes
sense since the degree of the line bundle is well defined on large
open subsets. This definition works for reductive groups as well.

\edefe

\bdefe\label{ramschm} A $H$--quasibundle $\cp$ is said to be
{\it ${\mu}$--semistable} (resp. {\it ${\mu}$--stable}) if the induced
rational principal $H$--bundle ${\cp}|_{U(\cp)}$ is {\it
${\mu}$--semistable}(resp.{\it ${\mu}$--stable}).\edefe

\bdefe\label{admis} A reduction of structure group of $E$ to a
parabolic subgroup $P$ is called {\it admissible} if for any character
$\chi$ on $P$ which is trivial on the center of $H$, the line bundle
associated to the $P$-bundle $E_P$ obtained by the reduction of
structure group, has degree zero.  \edefe

\bdefe\label{poly} An $H$-bundle $E$ is said to be {\it polystable} if
it has a reduction of structure group to a Levi subgroup $R$ of a
parabolic $P$ such that the $R$-bundle $E_R$ obtained by the
reduction, is stable and the extended $P$ bundle $E_R(P)$ is an
admissible reduction of structure group for $E$. Since the definition
involves only degrees of line bundles, it clearly holds good for
rational principal bundles as well. (cf. \cite{rr}) \edefe

\brem\label{quasipoly1} It is clear that we have the natural notions
of polystability for $H$--quasibundles as well. So a quasibundle
${\cp}$ is termed polystable if the induced rational principal
$H$--bundle ${\cp}|_{U(\cp)}$ is {\it polystable}.\erem

\subsubsection{Semistability and polystability over curves}

The study of semistability and polystability of principal bundles on
curves was initiated by A.Ramanathan. Over the years it has developed
in many directions and the results which one requires are scattered in
the literature. Polystability is also differently called {\it
  quasi-stability} in the literature but we avoid this terminology for
the obvious reasons (we have already a notion of quasibundles ..).

For the convenience of the reader we gather
some of the relevant facts with appropriate references in the
following theorem:

\bth\label{majorrecap} Let $C$ be a smooth projective curve over the
field of characteristic zero. Let $H$ a semisimple algebraic group.

The following are equivalent:
\begin{enumerate} 

\item[(i)] A principal $H$--bundle $E$ is polystable in the sense of
Ramanathan (Def \ref{poly}).
\item[(ii)] There exists a faithful representation $H \hra GL(V)$ such that
the associated bundle $E(V)$ is polystable of degree
$0$.
\item[(iii)] For every representation $H \lr GL(W)$, the bundle $E(W)$ is
polystable of degree $0$.
\item[(iv)] Let $ad :H \lr GL({\goth H})$ be the adjoint
representation. Then $E({\goth H})$ is polystable (of degree zero).
\item[(v)] The bundle $E$ arises from a representation $\rho: {\pi}_1(C)
\lr K$, where $K$ is a maximal compact subgroup of $H$.
\item[(vi)] The bundle $E$ carries an Einstein--Hermitian connection.
\end{enumerate}
\eeth

{\it Proof:} The equivalence of (i) and (v) is the main theorem of
Ramanathan, which generalises the Narasimhan--Seshadri theorem for
principal bundles. (\cite{rama}).

(i) $\Leftrightarrow$ (vi) is the main theorem of \cite{rs}. The
equivalence (iv) $\Leftrightarrow$ (v) is shown in \cite[Lemma
10.12]{atiyahbott}.

By \cite[Prop 1]{rs}, we may go modulo the center and assume that the
group $H$ has trivial center. Therefore, the adjoint representation
$ad : H \hra GL({\goth H})$ is a faithful representation. From this
standpoint, (ii), (iii) and (iv) are equivalent by a Tannakian
principle and the proof can be found in \cite[Prop 2.3]{basa}. The
argument there is for semistability, but the changes needed to be made
for polystability are easy since all bundles involved are of degree
$0$(see also the proof of Prop \ref{ssnonred} below).
\begin{flushright}{\it q.e.d}\end{flushright}

\brem\label{atibo} In fact, the equivalence [iv] $\Leftrightarrow$ [v]
holds more generally even for $H$ reductive (cf. \cite[Lemma
10.12]{atiyahbott}). \erem

\subsection{\it Some key lemmas}

We recall the following couple of facts about torsion--free sheaves
which we will use in this work. The first one is rather well known.

Let $C \subset U$ be a smooth projective curve. We recall:

\blem\label{dm1}(cf. \cite [Lemma 2.10]{basa} ) Let $T = \spec A$ and
let $E_T$ be a family of vector bundles on $C \times T$ such that
$E_p$ is semistable of degree $0$. Let $s_K$ be a section of the
family $E_K$ restricted to $C \times \spec K$, with the property that
for a base point $x \in C$, the section $s_K$ extends along $x \times
T$ to give a section of $E_{x \times T}$. Then the section extends to
the whole of $C \times T$.  \elem 

We have more generally:

\blem\label{dm2} Let $W$ be a family of semistable vector bundles with
$c_1 = 0$ on $U$ parametrised by $T$, i.e a vector bundle on $U_A$,
where $U \subset X$ is a {\it large} open subset of $X$. Suppose that
we are given a regular section $s_K : U_K \lr W_K$ such that, for an
irreducible smooth divisor $Y \subset U$, the section extends as a
regular section along $Y_A$. Then the section $s_K$ extends to a
regular section $s_A$ on $U_A$.  \elem

{\it Proof}: For the section $s_{T-p}$, viewed as a section of
$W_{T-p}$ we have two possibilities, since $U_A$ is normal and the
polar set is a divisor:
\begin{enumerate}
\item[(a)] it either extends as a regular section $s_T$.
\item[(b)] or it has a pole along $U \times p$.
\end{enumerate}

By the given property, we have a section $\sigma_{_Y}: Y \times T \lr
W_Y$ with the property that, at $\forall x \in Y$, $s_t (x) =
\sigma_{_Y}(x,t)$ $\forall ~ t \in T-p$;

So to complete the proof, we need to check that the
possibility (b) cannot hold:

Suppose it does hold. Then the section $s_{T-p} = s_K$ is a section of
$W_K$, i.e, a {\it rational section} of $W$ with a pole along the
divisor $U \times p \subset U \times T$, of order $k \geq 1$.

Thus by multiplying $s_{T-p}$ by $\pi^k$ we get a regular
section $s_T'$ of $W$ on $U \times T$.  If $s_T' = \{
s_t' \}_{t \in T}$, then we have:

{\renewcommand{\labelenumi}{{\rm (\roman{enumi})}}
\begin{enumerate}
\item $s_t' = \lambda (t) \cdot s_t$, $t \in T-p$ where
$\lambda : T \lr {\mathbb C}$ is a function given by $\pi^k$,
having zeros of order $k$ at $p$.
\item $s'_p$ is a non-zero section of $W_p$. Here we first note that
$W_p$ is a bundle on $U$ which is large and hence the
$\mu$--semistability is completely determined by $W_p$. Further, by
taking $W_p^{**}$, this extends as a reflexive sheaf to the whole
variety $X$ and remains $\mu$--semistable of degree $0$. By Hartogs'
theorem (since reflexive sheaves are {\it normal}), the section $s'_p$
extends to a non-zero section of $W_p^{**}$ to the whole of $X$.
\end{enumerate}}
\noindent
Since ${\mathcal O}_X$ is stable and degree $0$, it follows that
$s'_p$ gives a short exact sequence of sheaves on $X$:
\[
0 \lr {\mathcal O}_X \lr W_p^{**} \lr Q \lr 0
\]
where the torsion part of $Q$ is supported in a subset of codimension
$\geq 2$ (if the quotient sheaf had a divisor in its support then by
taking determinants, we see that there is a contradiction to the
equality of the degrees $deg({\mathcal O}_X) = deg(W_p^{**}) = 0$).

Thus it follows that there exists 
\[
y \in Y~ such~that ~s'_p(y) \neq 0. \eqno(*)
\]
\noindent
By the assumption that the section extends along $Y_A$, and by the
fact that $y \in Y$ it follows that $s_t (y) = \sigma_{_Y}(y,t)$, for $t
\in T-p$, and hence
$$s'_t (y) = \lambda (t) \cdot \sigma_{_Y}(y,t)~for~ t \in T-p.$$
Therefore, by continuity, since $\sigma_{_Y}(y,p)$ is well-defined,
we see that
$\lambda (t) \cdot \sigma_{_Y}(y,t)$ tends to $\lambda (p) \cdot
\sigma_{_Y}(y,p) =0$, as $t \lr p$.

Also, $s'_t (y) \lr s'_p (y)$ since $t \lr p$.  Hence by continuity
again, it follows that $s'_p (y) =0$, which contradicts $(*)$.

Thus the possibility (b) does not occur and we are
done.\begin{flushright} {\it q.e.d} \end{flushright}

\subsubsection{A Refinement of Langton's theorem}
Let $R^{\mu ss}$ denote the subset of the usual Quot scheme which
parametrises the usual torsion--free sheaves with fixed topological
type and Hilbert polynomial.  Let $\cg$ be the group $SL(\chi)$
($\chi$ being the Euler characteristics of sheaves in $R^{\mu ss}$)
acting on the (open subset)$R^{\mu ss}$ of the Quot scheme.

\blem\label{langt}{\it (Modified Langton)} Let $E_K$ denote a family
of $\mu$--semistable torsion--free sheaves of degree zero on $X \times
\spec K$. Then, (by going to a finite cover $S$ of $T$ if need be )
the sheaf $E_K$ extends (upto isomorphism) to a family $E_A$ with the
property that the limiting sheaf $E_0$ is in fact {\it polystable}.
\elem

{\it Proof.} By definition $E_K$ gives a $K$--valued point $x_K :
Spec(K) \lr R^{\mu ss}$

By Langton, there exists an $A$--valued point $x_A : Spec(A) \lr
R^{\mu ss}$ such that $x_p$ is given by $E_p$, which is a semistable
torsion--free sheaf.

One knows that (for a choice of a Jordan--Holder filtration) there
exists a family $F_t|_{t \in \ba}$ such that $F_t \simeq E_p$ for $t
\neq 0$ and the limit is a polystable sheaf $F_0 \simeq
gr^{\mu}(E_p)$. Let us denote this family by a morphism $f: {\ba} \lr
R^{\mu ss}$ and the point $f(0) = y_0$ for $t = 0$.

By going to a finite cover if need be, we may assume that there exists
a point $g \in {\cg}(\bc(t))$ such that $g \cdot E_p \simeq
F_{\bc(t)}$ with limit given by $F_0$.

Let $D = {\overline {{\cg} \cdot x_K}}$ (the ${\cg}(K)$--orbit closure
in $R^{\mu ss}$) considered as a $\bc$--variety. Then by definition
the point $x_p \in D$ (where $x_p$ corresponds to the sheaf
$E_p$). Further, $D$ is $\cg$--stable. Therefore, $g \cdot x_p \in
D({\bc}(t))$. Since $D$ is closed it follows that the point $y_0$
belongs to $D$.

Again, since $D$ is irreducible, we can join $y_0$ to the orbit ${\cg}
\cdot x_K$ and we get a scheme $S$, the spectrum of a discrete
valuation ring, such that $s_0 = y_0$ and $S - s_0 = Spec(L)$ where $L$
is a finite extension of $K$ (consider the natural map $q: {\cg} \cdot
x_K \lr Spec(K)$. Therefore $(S \cap {\cg} \cdot x_K )$ is mapped to
$Spec(K)$ by $q$ which gives the extension $L/K$.)

We also get a resulting family $\{E'_s\}_{s \in S}$ such that $E'_s
\simeq E_{q(s)}$ for $s \neq 0$ and $E'_{s_0} \simeq gr^{\mu}(E_p)$. We
are done. \begin{flushright} {\it q.e.d} \end{flushright}

\section{Extension of structure group to the flat closure}
 Fix a faithful representation $\rho: H \hookrightarrow G$ defined
 over ${\bc}$. Consider the extension of structure group of the
 principal $H_K$--quasibundle $(E_K, \sigma_K)$ via the induced
 $K$-inclusion $\rho_K: H_K \hookrightarrow G_K$. In other words, on
 a large open $U_K$, we are given a reduction of structure group of
 the principal $G_K$--bundle $E_K|_{U_K}$ to $H_K$.

Then, since $G = GL(n)$, by the properness of the functor of
 semistable torsion--free sheaves (the Theorem of Langton), there
 exists a {\it semistable extension} of $E_K$ to a torsion--free sheaf
 on $X \times \spec A$, which we denote by $E_A$.  Call the
 restriction of $E_A$ to $X \times p$ (identified with $X$) the {\it
 limiting bundle} of $E_A$ and denote it by $E_p$ (as in \S1). One has
 in fact slightly more, which is what we need.

\subsection{\it Mehta--Ramanathan restriction theorems}

Fix an integer $c \in [1,d-1]$ where $d = dim(X)$ and
\[
(a_1, \ldots, a_c) \in {\bz}^c
\]
with $a_i \geq 2$ for all $i \in [1,c]$. For any integer $m \geq 1$,
by a complete intersection of type $(m)$ we mean a complete
intersection of divisors:
\[
D_1 \cap D_2 \cap \cdots \cap D_c
\]
with $D_i \in |a_i^{m} \Theta|$. 

By a general complete intersection subvariety of type $(m)$ we mean
the complete intersection formed by the $D_j's$ from a non--empty open
subset $S$ of the linear system $\prod_{i = 1}^{n-1} |{a_i^{m} \Theta|}$.

Note that for $m \gg 0$, a {\it general complete intersection of type
$(m)$} is a smooth irreducible subvariety of $X$ of dimension $d - c$.

We recall the main theorem of \cite{mr}:

\bth (Mehta--Ramanathan). Let $V$ be a semi--stable
(resp. stable, polystable) rational vector bundle over $U \subset
X$. Then for all $m \gg 0$ the restriction $V|_{Y^o}$ to $Y \cap U$,
for a general complete intersection $Y$ of type $(m)$ is a
semi--stable (resp. stable, polystable) rational vector
bundle. Conversely, if $V|_{C}$ to a general complete intersection
curve of type $(m)$ is semi--stable (resp. stable, polystable) then so
is $V$. \eeth

\brem\label{crucial} Consider the extended sheaf $E_A$ obtained above
with the property that $E_p$ is a polystable torsion--free sheaf of
degree $0$. Since $E_p$ is a {\it polystable} torsion--free sheaf of
degree $0$, let $U = U(E_p) \subset X$ be the largest open set where
it is locally free. By the Restriction theorem of Mehta-Ramanathan,
for a large $m$, there is an open subset $S \subset \prod_{i =
1}^{n-1} |{a_i^{m} \Theta|}$, such that if $Y \in S$ is a a complete
intersection subvariety, then $E_p|_{Y \cap U}$ remains polystable of
degree $0$.

Let the $'$ denote intersection with $U$. We can therefore chose a
chain of subvarieties:
\[
C \subset \cdots (D_1 \cap D_2 \cap \cdots \cap D_{c})' \cdots \subset D_1'
\subset U \eqno(\ref {crucial})
\]
and the point $x \in C$ such that the restriction $E_p|_{(D_1 \cap D_2
\cap \cdots \cap D_{c})'}$ for every $c$ is {\it locally free and
remains polystable} of degree $0$.

Since $A$ is a discrete valuation ring, if we choose an open subset $U
\subset X$ where the limiting sheaf $E_p$ is locally free, then the
family $E_A$ is also {\it locally free} when restricted to the large
open subscheme $U_A$ (this is easy to see. cf. for example Lemma 5.4
\cite{newstead}). Hence, if we further restrict $E_A$ to the curve
$C_A \subset U_A$, we get a family of locally free sheaves $E_A
|_{C_A}$ on the smooth projective curve $C$ parametrised by $Spec
A$. Further, by choice, the limiting bundle $E_p$ is polystable on $C$
of degree $0$. \erem

\subsection{\it The Flat closure}
We observe the following:
\begin{itemize}

\item Note that giving a reduction of structure group of the
$G_K$-bundle $E_K$ on a large $U_K$ is equivalent to giving a section
$s_K$ of $E_K (G_K/H_K)$ over $U_K$ .

\item Let
\[
C \subset \cdots (D_1 \cap D_2 \cap \cdots \cap D_{c})' \cdots \subset
D_1' \subset U = U(E_p) 
\]
be as in ($\ref{crucial}$) above. We fix a base point $x \in C
\subset U$ and denote by $x_A = x \times \spec A$, the induced section
of the family (which we call the {\it base section}):
\[
X_A \lr \spec A
\]
Since $x \in U$, it allows us to work with genuine principal bundles
and their restrictions to the section $x_A$ defined by the base point
$x$.

\item {\sf From now on, unless otherwise stated, we shall fix this chain of
smooth subvarieties of $U$ along with the base point $x$.}

\item Let $E_{x,A}$ (resp $E_{x,K}$) be as in \S1, the restriction of
$E_A$ to $x_A$ (resp $x_K$). Thus $s_K(x)$ is a section of
$E_K(G_K/H_K)_x$ which we denote by $E_x (G_K/H_K)$.

\item Since $E_{x,A}$ is a principal $G$-bundle on $\spec A$ and
therefore trivial, it can be identified with the group scheme $G_A$
itself. {\it For the rest of the article we fix one such
identification, namely:}
\[
\xi_A : E_{x,A} \longrightarrow G_A.
\]
\item Since we have fixed $\xi_A$ we have a canonical identification
\[
E_x(G_K/H_K) \simeq G_K/H_K 
\]
which therefore carries a natural {\it identity section} $e_K$ (i.e
the coset $id.H_K$). Using this identification we can view $s_K(x)$ as
an element in the homogeneous space $G_K/H_K$.

\item Let $\theta_K \in G(K)$ be such that $\theta_K^{-1} \cdot s_K(x)
  = e_K$.  Then we observe that, the isotropy subgroup scheme in $G_K$
  of the section $s_K(x)$ is $\theta_K.H_K.\theta_K^{-1}$. (We remark
  that such a $\theta_K$ will exist after going to a finite extension
  of $K$. By an abuse of notation we will continue calling this
  extension as $K$. This is required since $G_K \longrightarrow
  G_K/H_K$ need not be locally trivial ).

\item On the other hand one can realise $s_K(x)$ as the identity coset of
$\theta_K.H_K.\theta_K^{-1}$ by using the following
identification:
 $$G_K/\theta_K.H_K.\theta_K^{-1} 
 \stackrel{\sim}{\longrightarrow} G_K/H_K.$$
 $$g_K(\theta_K.H_K.\theta_K^{-1}) \longmapsto
 g_K\theta_K.H_K.$$ 
\end{itemize}
\noindent 
\bdefe Let ${H_K'}$ be the subgroup scheme of $G_K$ defined as: 
\[
H_K' := \theta_K.H_K.\theta_K^{-1}.
\]
\edefe
\noindent
Using $\xi_A$ we can have a canonical identification:
\[
 E_{x}(G_K/H_K') \simeq G_K/H_K'.
\] 
\noindent
Then we observe that, using the above identification we get a section
$s_K'$ of $E_K (G_K/H_K')$, with the property that, $s_K'(x)$ is the
{\it identity section} and moreover, since we have conjugated by an
element $\theta_K \in G_A(K) (=G(K))$, the isomorphism class of the
$H_K$-bundle $P_K$ given by $s_K$ does not change by going to $s_K'$.

Thus in conclusion, the $G_A$-bundle $E_A$ has a reduction to
$H_K'$ given by a section $s_K'$ of $E_K (G_K/H_K')$, with the
property that, at the given base section $x_A = x \times \spec A$, we
have an equality $s_K' (x_A) = e_K'$, with the {\it identity element}
of $G_K/H_K'$ (namely the coset $id.H_K'$).

\bdefe\label{fc} The {\it flat closure} of the reduced group scheme
$H_K'$ in $G_A$ is defined to be the schematic closure of $H_K'$ in
$G_A$ with the reduced scheme structure.  Let $H_A'$ denote the
flat-closure of $H_K'$ in $G_A$.
 \edefe

\brem\label{onflatclosure} Let $I_K$ be the ideal defining the
 subgroup scheme $H_K'$ in $K(G)$ (note that $G_A$ (resp $G_K$) is an
 affine group scheme and we denote by $A(G)$ (resp $K(G)$) its
 coordinate ring). If we set $I_A = I_K \cap A(G)$, then it is easy to
 see that since we are over a discrete valuation ring, $I_A$ is in
 fact the ideal in $A(G)$ defining the flat closure $H_A'$. \erem

\noindent
We then have a canonical identification via $\xi_A$:
\[
 E_{x}(G_A/H_A') \simeq G_A/H_A'.
\] 
\noindent
One can easily check that $H_A'$ is indeed a subgroup scheme of $G_A$
since it contains the identity section of $G_A$, and moreover, it is
faithfully flat over $A$. Notice however that $H_A'$ {\it need not} be
a {\it reductive} group scheme; that is, the special fibre $H_k$ over
the closed point need not be reductive.

Observe further that $s_K'(x)$ extends in a trivial fashion to a
 section $s_A'(x)$, namely the {\it identity coset section} $e_A'$ of
 $E_x(G_A/H_A')$ identified with $G_A/H_A'$ .

\subsection{\it Chevalley embedding of $G_A/H_A'$}

 As we have noted, $H_A'$ need not be reductive and the rest of the
 proof is to get around this difficulty.  Our first aim is
 to prove that the structure group of the bundle $E_A(G_A)$ can be
 reduced to $H_A'$ which is the statement of Proposition
 \ref{flatext}.

 We need the following generalisation of a well-known result of
 Chevalley from \cite{basa}.

\blem\label{chev} 
There exists a finite dimensional $G_A$-module $W_A$ such that
 $G_A/H_A'~\hookrightarrow~W_A$ is a $G_A$-immersion.  
\elem
 
\subsection{\it Extension to flat closure and local constancy}

Recall that the section $s_K'(x)$ extends along the base section
$x_A$, to give $s_A'(x) = w_A$. The aim of this section is to prove the
following key theorem.
\noindent 
\bth\label{flatext} There exists a {\it large} open subset $U^o
 \subset U$ such that the section $s_K'$ extends to a section $s_A'$
 of $E_A(G_A/H_A')$ when we restrict $E_A$ to the open subset $U^o_A
 \subset U_A$. In other words, the structure group of
 $E_A{_{|_{U^o_A}}}$ can be reduced to $H_A'$; in particular, if $H_k'$
 denotes the closed fibre of $H_A'$, then the structure group of $E_k$
 can be reduced to $H_k'$.  \eeth
 
 Towards this we need the following key result.

\bprop\label{lc} Let $E$ be a {\it polystable} principal $G$-bundle of
 degree zero on a smooth projective curve $C$ (here $G$~=~$SL(n)$~ or
 $GL(n)$).  Let $W$ be a $G$-module and $N$ a $G$-subscheme of $W$ of
 the form $G/H'$.

 If $s$ is a section of $E(W)$ such that for some $x \in C$, $s(x)$
 lies in the fibre $E(N)_x$ of the fibration $E(N)$ at $x \in C$, then
 the entire image of $s$ lies in $E(N)$. Consequently, $E$ has a
 reduction of structure group to the subgroup $H'$.

\eprop

{\it Proof.} The bundle $E$ being polystable, it is defined by a
 ``unitary'' representation
\[
\chi : \pi_1(C) \longrightarrow G
\]
which maps into a maximal compact subgroup of $G$. This implies that
if the universal covering $j : Z \longrightarrow C$ is considered as a
principal fibre space with structure group $\pi_1(C)$, then the
principal $G$-bundle $E$ is the associated bundle through $\chi$.

Let $\rho : G \longrightarrow GL(W)$ be the representation defining
the $G$-module $W$. Then $E(W)$ can be considered as the bundle
associated to the principal bundle $j : Z \longrightarrow C$ through
the representation
\[
\rho \circ \chi :\pi_1(C) \longrightarrow GL(W).
\]
which maps into the unitary subgroup of $GL(W)$. 

By generalities on principal bundles and associated constructions,
since 
\[
E(W) \simeq Z \times^{\pi_1(C)} W
\]
a section of $E(W)$ can be viewed as a $\pi_1(C)$-map
\[
s_1 : Z \longrightarrow W
\]
Now, since $Z$ is the universal cover of the curve $C$ and $s$ is a
section of $E(W)$, then one knows (cf. \cite{ns}) that there exists a
$\pi_1(C)$-invariant element $w \in W$ such that $s$ is defined by a
map
\[
s_1 : Z \longrightarrow W
\]
given by $s_1(x) = w ,~ \forall x \in C$, i.e ``the constant map
sending everything to $w$''.

Since $w \in W$ is a $\pi_1(C)$-invariant vector and the action of
$\pi_1(C)$ is via the representation $\chi$, we see that $\chi$
factors via
\[
\chi_1: \pi_1(C) \lr H'
\]
since $H' = Stab_G(w)$.

In particular, we get the $H'$-bundle from the representation $\chi_1$
and clearly this $H'$-bundle is the reduction of structure group of
the $G$-bundle $E$ given by the section $s$.

By the very construction of the reduction, the induced $H'$-bundle is
{\it flat} and also semistable since it comes as the reduction of
structure group of the polystable bundle $E$ (by Def \ref{ssnonred1}
below). This proves the Proposition \ref{lc}. \begin{flushright} {\it
q.e.d} \end{flushright}

\bdefe\label{ssnonred1} Let $H'$ be an affine algebraic group not
necessarily reductive. Let $P$ be a principal $H'$-bundle on $C$ a
smooth projective curve. We define $P$ to be {\it polystable}) if: 

\begin{enumerate}
\item It is {\it flat} (in the sense that it comes from the
representation $\chi$ of the fundamental group of $C$)

\item there exists a faithful representation
\[
\rho :H' \lr GL(V)
\]
such that the associated vector bundle $P(V)$ is {\bf polystable of
degree zero}).  
\end{enumerate}
\edefe

\brem This definition is ad hoc and made only to suit our needs.  Let
${\cm}:= {\cm}({\rho} \circ {\chi})$ denote the Zariski closure of the
image of ${\rho} \circ {\chi}: {\pi}_1(C) \lr GL(V)$ in $GL(V)$.  We
term this the ``monodromy'' subgroup associated to the representation
${\rho} \circ {\chi}$.

Since the bundle associated to ${\rho} \circ {\chi}: {\pi}_1(C) \lr
GL(V)$ i.e $P(V)$, is assumed to be polystable of degree $0$, by the
Narasimhan--Seshadri theorem, the representation ${\rho} \circ {\chi}$
is unitary and also the monodromy subgroup ${\cm}$ is {\it reductive}
(possibly non-connected). This can be viewed as a Tannakian
interpretation of polystability. Further, since $\chi$ maps into $H'$,
the monodromy subgroup $\cm$ is a subgroup of $H'$.

Let us denote the inclusion of the monodromy subgroup by:
\[
{\iota}: {\cm} \hra H'
\]
The fact that the bundle $P$ comes from an associated construction via
the homomorphism $\chi$ implies that $P$ has a reduction of structure
group to ${\cm}$. Let the resulting ${\cm}$ bundle be denoted by
$P_{_{\cm}}$ 

\erem

\bprop\label{ssnonred} Let $H'$ be an affine algebraic group (not
necessarily reductive), as above, and let $P$ be a {\it polystable}
principal $H'$-bundle on a smooth projective curve $C$. Let $f :H' \lr
H$ be a morphism from $H'$ to a semisimple group $H$. Then the
associated principal $H$-bundle $P(H)$ is also polystable.  \eprop

{\it Proof.} Since $P$ has a reduction of structure group to ${\cm}$,
we may view the principal $H$ bundle $P(H)$ as obtained from the
homomorphism $f \circ {\iota}:{\cm} \lr H$. Thus $P(H)$ can be
identified with $P_{_{\cm}}(H)$. 

To check the polystability of the principal $H$--bundle $P(H)$ by
Theorem \ref{majorrecap}, we need only check that if $\psi :H \lr
GL(W)$ is any representation of $H$ then the associated bundle
$P(H)(W)$ is polystable. In other words we need to check that if
$\gamma :{\cm} \lr GL(W)$ (with $\gamma = \psi \circ f \circ {\iota}$)
then the associated bundle $P_{_{\cm}}(W) = P(H)(W)$ is polystable.

Observe that by Def \ref{ssnonred1} we have a faithful representation
$GL(V)$ of ${\cm}$ (from that of $H'$) such that $P_{_{\cm}}(V)$ is
polystable of degree zero. Further, we are over a field of
characteristic zero and so the ${\cm}$-module $W$ can be realised as a
direct summand of a direct sums of some $T^{a,b}(V) = V^{\otimes a}
\otimes {V^{*}}^{\otimes b}$ (cf. for e.g \cite{sim} pp 86. We need the
reductivity of ${\cm}$ here, otherwise in general this is only a
subquotient and not a direct summand).

Hence the vector bundle $P_{_{\cm}}(W)$ is a direct summand of
$\bigoplus T^{a,b}(P_{_{\cm}}(V))$. Since $P_{_{\cm}}(V)$ is
polystable of degree zero, so is $T^{a,b}(P_{_{\cm}}(V))$.

By assumption, since $P_{_{\cm}}$ is {\it flat} the associated vector
bundle via any representation is of {\it degree zero}. Hence
$P_{_{\cm}}(W)$ has degree $0$.

Since all bundles have degree zero and $\bigoplus
T^{a,b}(P_{_{\cm}}(V))$ is a direct sum of stable bundles of degree
$0$ it is easy to check that $P_{_{\cm}}(W)$ is also direct sum of
stable bundles of degree $0$ and hence polystable.\begin{flushright}
{\it q.e.d}\end{flushright}

\brem The polystability of $P_{_{\cm}}(W)$ also follows directly from
that of $P_{_{\cm}}$ by Remark \ref{atibo}. We cannot use the Theorem
\ref{majorrecap} as it stands since $\cm$ is in general only
reductive. But the advantage here is that the bundle $P_{_{\cm}}$
comes from a representation of the fundamental group.\erem

After this brief interlude on curves, we now return to the general
setting of higher dimensional varieties $X$.

\bprop\label{lcs} Assume that $dim(X) = 2$. Let $E$ be a rational {\it
 polystable} principal $G$-bundle of degree zero on $U \subset X$.
 Let $W$ be a $G$-module and $N$ a $G$-subscheme of $W$ of the form
 $G/H'$ where $H' = Stab_G(w)$ for some $w \in W$. If $s$ is a section
 of $E(W)$ such that for $x \in C \subset U$, the image $s(x) = w$ in
 the fibre of $E(N)$ at $x \in X$. Then there exists a {\it big} open
 subset $U^o \subset U$, such that the entire image of $U^o$ under $s$
 lies in $E(N)$. (Here $C \subset U$ is as in $({\ref{crucial}})$
 where the base point $x$ was chosen in a general $C$). In particular,
 $E$ has a reduction of structure group to the subgroup $H'$ on a big
 open subset.

\eprop

{\it Proof}: By choice, since $x \in C \subset U$, and $E|_C$ remains
polystable of degree $0$, it follows by Prop \ref{lc} that 
$s(C) \subset N$.

Now, again by the Mehta-Ramanathan theorem for $m \gg 0$, there exists
an open $\Omega_m \subset |m {\Theta}|$ such that for $t \in \Omega_m$, the
restriction of $E|_{Y_t}$ to the smooth projective curve $Y_t$ remains
{\it polystable} of degree $0$. Further, the set $\{t \in \Omega_m~|~ Y_t
\cap C \neq \emptyset \} $ is an open subset of $\Omega_m$, which we
continue to call $\Omega_m$. Since $Y_t \cap C \neq \emptyset $, for each
$t \in \Omega_m$, there is a point $x_t \in Y_t \cap C$ and hence $s(x_t)
\in Y_t$. Thus applying Prop \ref{lc} to the restriction $E|_{Y_t}$,
we see that $s(Y_t) \subset N$ for all $t \in \Omega_m$. Hence, the dense
subset $U_1 = \bigcup Y_t \subset U$ is mapped into $N$ i.e $s(U_1)
\subset N$. This implies that the entire image $s(U)$ lies in the
closure of the orbit ${\overline N} \subset W$.

Now observe that $N \subset {\overline N}$ is an open subset,
therefore, it follows that the set of points $U^o = \{u \in U | s(u)
\in N \}$ is an open subset of $U$.

We claim that $U^o \subset U$ is also {\it large}. For if not, then
the complement of $U_2$ in $U$ contains a curve and there are $t \in
\Omega_m$ such that this curve in the complement meets $Y_t$, a
contradiction since the entire curve $Y_t$ gets mapped into $N$. This
proves our claim. \begin{flushright} {\it q.e.d} \end{flushright}

Now let $X$ be an smooth projective variety of dimension $d$.

Since $E$ is a rational polystable bundle there is a {\it big} $U$
where it is a principal $H$--bundle. We work with the chosen complete
intersection as in $(\ref{crucial})$ in such a manner that $C \subset
U$. With these choices and the point $x \in C \subset U$ we have the
following:

\bprop\label{lcs1} Let $E$ be a rational {\it polystable} principal
 $G$-bundle of degree zero on $U \subset X$.  Let $W$ be a $G$-module
 and $N$ a $G$-subscheme of $W$ of the form $G/H'$ where $H' =
 Stab_G(w)$ for some $w \in W$. Let $s$ is a section of $E(W)$
 (defined on $U$) such that the image $s(x)$, of the point $x \in C
 \subset U$, lies in the fibre of $E(N)$ at $x \in X$.  Then there
 exists a {\it big} open subset $U^o \subset U$, such that the entire
 image of $U^o$ under $s$ lies in $E(N)$. In particular, $E$ has a
 reduction of structure group to the subgroup $H'$ on a big open
 subset.

\eprop

{\it Proof}: By choice, since $x \in C \subset U$, and $E|_C$ remains
polystable of degree $0$, it follows by Prop \ref{lc} that $s(C)
\subset E(N)$.

Let us, for the present, denote the divisor $D_1$ simply by $D$.  By
an induction on the dimension $d$, and Prop \ref{lcs}, we see that
since $C \subset D' = D \cap U$, and since $E|_{D'}$ is also
polystable, there exists a big open subset $D^o \subset D' $ such that
$s(D^o) \subset E(N)$.

Let $Y_j \in S$ be any other general divisor. Then, by the property of
$S$, we see that $D \cap Y_j$ is a smooth divisor of $D$ and
$Y_j$. Further, by other choices of $Y_i \in S$ we see that there
exists a smooth projective curve $C_j \subset Y_j \cap D$. 

Since $D^o \subset D'$ is a {\it big} open subset, in fact, it is easy
to see that, by going to a smaller open subset of $S$ if need be (and
fixing it), we can make sure that $C_j \subset D^o$ for all $Y_j \in
S$ (cf. 4.4, \cite{mr}).

Again $C_j \subset (Y_j \cap U)$ and $E|_{(Y_j \cap U)}$ is also a
polystable rational bundle. Since $s(C_j) \subset s(D^o) \subset
E(N)$, it follows again by induction on dimension that there exists a
big open $Y_j^{o} \subset Y_j$ such that $s(Y_j^{o}) \subset E(N)$ for
all $Y_j \in S$.

Now as in Prop \ref{lcs}, we see that the set of points $U^o = \{u \in
U | s(u) \in E(N)_u \}$ is a non--empty open subset of $U$.

We claim that $U^o \subset U$ is also {\it large}. For if not, then
the complement of $U^o$ in $U$ contains a divisor $R \subset U$ and
there are $Y_t \in S$ such that $R \cap Y_t \neq \emptyset$.
Moreover, by Bertini's theorem one has $dim(R \cap Y_t) = dim R - 1 =
dim(Y_t) - 1$. This implies that $(R \cap Y_t) \not\subset (Y_t -
Y_t^{o})$, i.e $R \cap Y_t^{o} \neq \emptyset$. That is, there exists $y
\in R \cap Y_t^{o}$ such that $s(y) \in E(N)_y$ which contradicts the
fact that $R \subset U - U^o$.
\begin{flushright} {\it q.e.d} \end{flushright}

\brem We work with this $U^o$ from now on and since there is no
confusion we will call this $U$.\erem

\subsection{\it Completion of proof of Theorem \ref{flatext}}

\noindent
We work over the large open $U_A$ and all bundles in this proof are
are over $U_A$

By Lemma \ref{chev}, we have an immersion
\[
E_A(G_A/H_A') \subset E_A(W_A).
\]  

The given section $s_K'$ of $E_K(G_K/H_K')$ therefore gives a section
 $u_K$ of $E(W_K)$.  

Recall also that we have a chosen complete intersection 
\[
C \subset \cdots (D_1 \cap D_2 \cap \cdots \cap D_{c})' \cdots \subset D_1'
\subset U
\]
and the point $x \in C$
(cf. Remark \ref{crucial}).

Further, $u_K(x)$, the restriction of $u_K$ to $x
 \times T^*$, extends to give a section $u_A (x)$ of $E_x (W_A)$
 (restriction of $E_A(W_A)|_{C_A}$ to $x \times T$).

Thus by Lemma \ref{dm1}, and by the semistability of $E(W_A)_p$,
 on $C$ the section $u_K$ extends to give a section $u_A$ of
 $E(W_A)|_{C_A}$ over $C \times T$.

Again, by Lemma \ref{dm2}, we see, by an induction on dimension, that
the section $u_K$ which is defined on $U_K$ extends to section $u_A$
of $E_A(W_A)$ over the entire open subscheme $U_A$.

Now, to prove the Theorem \ref{flatext} , we need to make sure that:
$$ 
\begin{array}{l}
\mbox{The image of this extended section $u_A$ actually lands in
$E_A(G_A/H_A')$} \\
\end{array} . \eqno (*)
$$
This would then define $s_A'$.

To prove ($\ast$), it suffices to show that $u_A(U \times p)$ lies in
$E_A(G_A/H_A')_p$ (the restriction of $E_A(G_A/H_A')$ to $U \times
p$).
 
Observe that, $u_A(x \times p)$ lies in $E_A(G_A/H_A')_p$ since
$u_A(x) = s_A'(x) = w_A$.
 
Observe further that, if $E_p$ denotes the principal $G$-bundle on
$U$, which is the restriction of the $G_A$-bundle $E_A$ on $U \times
T$ to $U \times p$, then $E_A(W_A)_p = E_A(W_A)|{(U \times p)}$, and we
also have

\[
\begin{array}{ccc}
E_A(G_A/H_A')_p &
\stackrel{\simeq}{\longrightarrow} &
E_p(G_k/H_k') \\
\Big\downarrow & &
\Big\downarrow \\
E_A(W_A)_p   & \stackrel{\simeq}{\longrightarrow} &
E_p(W)
\end{array}
\]
and the vertical maps are inclusions:
\[
E_A(G_A/H_A')_p \hookrightarrow E_A(W_A)_p,~~and~~ E_p(G_k/H_k')
\hookrightarrow E_p(W)
\]
where $E_p(W) = E_p \times^{H_k'} W$ with fibre as the $G$-module $W =
W_A \otimes k$.  Note that $G/H_k'$ is a $G$-subscheme $Y$ of $W$.  

\noindent
Recall that $E_p$ is polystable of degree zero.  Then, from the
foregoing discussion, the assertion that $u_A(U \times p)$ lies in
$E_A(G_A/H_A')$, is a consequence of Proposition \ref{lcs1} applied to
$E_p$ (we might have to throw away a subset of $U$ of codimension
$\geq 2$ for this).
 
Thus we get a section $s_A'$ of $E_A(G_A/H_A')$ on $U \times T$, which
extends the section $s_K'$ of $E_A(G_A/H_A')$ on $U \times T^*$.
This gives a reduction of structure group of the $G_A$-bundle $E_A$ on
$U \times T$ to the subgroup scheme $H_A'$ and this extends the given
bundle $E_K$ to the subgroup scheme $H_A'$.
\begin{flushright} {\it q.e.d} \end{flushright}

\section{Semistable reduction for quasibundles over projective varieties}

The aim of this section is to prove the following theorem.
 
\bth\label{ssred} Let ${\cp}_K$ be a family of semistable principal
$H$-quasibundles on $X \times \spec K$, or equivalently, if $H_K$
denotes the group scheme $H \times \spec K$, a semi-stable
$H_K$-quasibundle ${\cp}_K$ on $X_K$. Then there exists a finite
extension $L/K$, with the integral closure $B$ of $A$ in $L$, such
that, ${\cp}_K$, after base change to $\spec B$, extends to a
semistable $H_B$-quasibundle ${\cp}_B$ on $X_B$.  \eeth
\noindent

\subsection{\it Potential good reduction}
 We begin by observing that we have extended the original rational
 $H_K$-bundle upto isomorphism to a rational $H_A'$-bundle. To
 complete the proof of the Theorem \ref{ssred}, we need to extend the
 $H_K$-bundle to an $H_A$-bundle.

\brem We note that in general the group scheme $H_A'$ obtained above
need not be a smooth group scheme over $A$. But in our case since the
characteristic of the base field is zero and since $H_A'$ is flat, it
is also smooth over $A$.  \erem

Recall that $H_A$ denotes the reductive group scheme $H \times \spec
A$ over $A$. We need the following crucial result from \cite{basa}. It
involves key inputs from Bruhat--Tits theory. For details see Appendix
of \cite{bapa}.

\bprop\label{serre} There exists a finite extension $L/K$ with the
following property: If $B$ is the integral closure of $A$ in $L$, and
if $H_B'$ are the pull-back group schemes, then we have a morphism of
$B$-group schemes
\[
H_B'~\longrightarrow H_B
\]
 which extends the isomorphism $H_L' \cong H_L$.
 \eprop

We also need the following:

\blem\label{ssnonred2} Let $H'$ be a non-reductive group and let $\psi
: H' \lr H$ be a homomorphism of algebraic groups. Let $P$ be a
rational principal $H'$--bundle on the big open subset $U \subset X$
obtained as reduction of structure group of a polystable vector bundle
$E$ of {\it degree zero} through an inclusion $H' \hra GL(V)$. Then,
by the extension of structure group the rational principal $H$--bundle
$P(H)$ is also semistable.\elem

{\it Proof:} By the choice of $P$ and by the restriction theorem of
Mehta and Ramanathan, we may choose a high degree general smooth
complete intersection curve $C$ of type $(m)$, with $C \subset U$ such
that $E|_{C}$ is polystable of degree zero. Let $a$ denote the product
$a_1 \ldots a_{d-1}$ in the notation of (3.1).

Thus by Prop \ref{ssnonred}, the associated $H$--bundle $P(H)|_{C}$ is
a semistable principal bundle on $C$.

We {\it claim} that this implies that $P(H)$ is itself a semistable
principal $H$--bundle. For, if $Q \subset H$ is a parabolic subgroup
and $\chi$ a dominant character of $Q$, $P(H)_{_Q}$ a $Q$--bundle
obtained from a reduction of structure group to $Q$, note that
\[
deg P(H)_{_Q}(\chi) \cdot a^m = deg P(H)_{_Q}(\chi)_{|_{C}}
\]
where $P(H)_{_Q}(\chi)$ denotes the line bundle associated to the
character $\chi$. Since $m \gg 0$ it follows by the semistability
of $P(H)_{|_{C}}$ that $deg P(H)_{_Q}(\chi)_{|_{C}} > 0$ and hence
$deg P(H)_{_Q}(\chi) > 0$, i.e $P(H)$ is
semistable. \begin{flushright} {\it q.e.d} \end{flushright}

{\it Proof.(of Theorem \ref{ssred}}) Let the semistable principal
$H_K$--quasibundle ${\cp}_K$ arise out of the faithful representation
$H \subset G = GL(V)$. In other words, there is a pair $({\ce}_K,
\sigma_K)$ of semistable torsion--free sheaf ${\ce}_K$ on $X_K$ and
section $\sigma_K$. By Prop \ref{langt}, we have a semistable
extension ${\ce}_A$ of ${\ce}_K$ (possibly by going to a finite
extension $L/K$ such that ${\ce}_p$ is polystable of degree $0$. Let
$U = U({\ce}_p)$. Then as we have seen earlier, the entire family
${\ce}_A$, when restricted to $U_A$ is locally free. Now we are in the
setting of Prop \ref{flatext}.

Moreover, the fibre of $E_B$ over the closed point is indeed {\bf
semistable}. To see this , observe firstly that it comes as the
extension of structure group of $E_p'$ by the map $\psi_k : H_k' \lr
H_k$. Recall (Proposition \ref{lc}) that $E_p'$ is the {\it
semistable} $H_k'$-bundle obtained as the reduction of structure group
of the polystable vector bundle $E(V_A)_p$ and the semistability of
the extended $H_k$--bundle follows by Lemma \ref{ssnonred2}. 

By Proposition \ref{flatext} we have a {\it rational} $H_A'$-bundle
$P_A'$ on $U_A$, which extends the $H_K$-bundle upto
isomorphism. Further, by Proposition \ref{serre}, by going to the
extension $L/K$ we have a morphism of $B$-group schemes $\psi_B : H_B'
\longrightarrow H_B$ which is an isomorphism over $L$. Therefore, one
can extend the structure group of the bundle $P_B'$ to obtain a
rational $H_B$-bundle $P_B$ which extends the $H_K$-bundle $P_K$.  By
using the faithful representation $H \subset GL(V)$, we get a family of
vector bundles $P_B(V)$ on $U_B$.

Note that the rational $H_K'$--bundle $P_K$ is the restriction to
$U_K$ of the principal $H_K$--quasibundle $({\ce}_K, \sigma_K)$.

In other words, we have an isomorphism $P_L(V) \simeq {\ce}_L$ on
$U_L$.

Now we apply a result from Langton's paper (\cite[Prop 6]{langton}) to
conclude that the bundle $P_B(V)$ on $U_B$ in fact extends to a
torsion--free sheaf ${\ce}_B$ on $X_B$. Further, the reduction section
$\sigma_L$ (defined on $X_L$ has extended to a section $\sigma_B$ on
$U_B$ (given by the principal $H_B$--bundle $P_B$). We can view the
reduction datum as follows:
\[
\sigma_B : (U_B \cup X_L) \lr (S_{{\ce}_B})//H_B
\]
(which by the extension property, agrees on the intersection). The
complement of $U_B \cup X_L$ in $X_B$ is a subset of points of
codimension $\geq 3$, and one can embed the affine scheme
$(S_{{\ce}_B})//H_B$ as a {\it closed subscheme} of a vector bundle
over $X_B$. 

This implies that the section $\sigma_B$ extends to a section
$\sigma_B : X_B \lr (S_{{\ce}_B})//H_B$. In other words, the pair
$({\ce}_K, \sigma_K)$ has been semistably extended (upto isomorphism)
to a pair $({\ce}_B, \sigma_B)$, i.e, the principal $H_K$--quasibundle
${\cp}_K$ has been semistably extended (upto isomorphism) to principal
$H_B$--quasibundle ${\cp}_B$, possibly after going to a finite
extension $B/A$. This completes the proof of the Theorem
\ref{ssred}. \begin{flushright} {\it q.e.d}
\end{flushright}

For purposes of applications later we isolate the following easier
half of the Mehta--Ramanathan restriction theorem: 

\blem\label{rest} Let $\cp$ be an $H$--quasibundle. Then $\cp$ is
$\mu$--polystable (resp stable) if ${\cp}|_{C}$ is so for a general
high degree smooth complete intersection curve $C \subset X$.  \elem

{\it Proof}: The proof is along the lines of the first part of the
proof of Prop \ref{ssnonred2}. For more general results along similar
lines cf.  \cite{bisgom}, whose methods work easily enough for
rational principal bundles and hence for quasibundles.

\section{Construction of the moduli space of bundles over surfaces}
From now onwards $X$ will be a {\it smooth projective surface}.  The
aim of this section is to give an algebro-geometric construction of
the Donaldson-Uhlenbeck compactification for the moduli space of
semistable principal $H$--bundles. This, in particular also describes
geometrically the Yang-Mills compactification of the moduli space of
ASD connections on principal bundles with arbitrary structure groups.
The method of proof is along the lines of the proof of J.Li (cf.
\cite{li}) and the methods in the paper of Le Potier (\cite{lepot})
(cf. also \cite[Chapter 8]{hl}).

\subsubsection{Double duals for quasibundles}
 If ${\cp}|_{U}$ is a rational $H$--bundle, where $U$ is a {\it big}
open subset of $X$, then by a theorem of Colliot-Th\'el\`ene and
Sansuc(\cite[Theorem 6.7]{collio}) there is a principal bundle (unique
upto isomorphism) which extends ${\cp}|_{U}$. This extension is
therefore also associated to each quasibundle $\cp$ as well.

\bdefe\label{doublestar}({\it Double duals}) We call this extended
principal bundle the double dual of the quasibundle $\cp$ and denote
it by ${\cp}^{**}$. This plays the role of the {\it double dual} of
torsion free sheaves. \edefe

\brem If $\cp$ is a $\mu$--semistable (resp. polystable,stable)
$H$--quasibundle, then since these definitions involve only {\it big}
open subsets, the corresponding double dual ${\cp}^{**}$ is a
$\mu$--semistable (resp. polystable,stable) principal $H$--bundle. \erem

\subsubsection{On determinant line bundles}

Recall that if $\cf$ is a flat family of coherent sheaves on $X$
parametrised by a scheme $S$, then $\cf$ defines an element $[\cf] \in
K^{0}( S \times X)$, the Grothendieck group of $S \times X$ generated
by locally free sheaves. We may then define the homomorphism from the
Grothendieck group of coherent sheaves on $X$ given by:
\[
{\lambda}_{\cf} : K(X) \lr Pic(S).
\]
This has a collection of functorial properties for which we refer to
(\cite{hl} page 179). For every class $u \in K(X)$ we denote the
associated class of the line bundle by ${\lambda}_{\cf}(u)$. Fix a
class $c \in K(X)$ with rank $r$ and Chern classes $c_1 = {\co}_X$ and
$c_2$, the very ample divisor $\Theta$ on $X$ and a base point $x \in
X$. Then there is a natural choice of a class $u_1(c) \in K(X)$
defined in terms of these fixed data (cf. page 184 \cite{hl}).

\subsubsection{The parameter space for quasibundles}

We first briefly recall the construction of the parameter space of
$\mu$--semistable torsion--free sheaves over the surface $X$. Let $c
\in K(X)$ be as above. Then one knows that the set $\cs$ of
isomorphism classes of $\mu$--semistable torsion--free sheaves of
class $c$ with trivial determinant and fixed Hilbert polynomial $P$ is
bounded and hence for a $m \gg 0$ we can realise them as points of a a
suitable quot scheme. Fix such an $m$. Let $W$ be a complex vector
space of dimension $P(m)$.

Let $R^{{\mu} ss} \subset Quot(W,P)$ be the locally closed subscheme
of all quotients $[q: W \otimes {\co}_X(-m) \lr {\ce}]$ such that
$\ce$ is $\mu$--semistable of rank $r$ with trivial determinant and
second Chern class $c_2$, with Hilbert polynomial $P$ and such that
$q$ induces an isomorphism $M \simeq H^0({\ce}(m))$.

There is a natural action of ${\cg} = SL(M)$ on $R^{{\mu}ss}$ and we
have the universal quotient $q_Q : {\pi}_X^{*}{\co}_{X}(-m) \otimes W
\lr {\ce}_Q$.

Returning to our setting, let $\rho: H \hra G$ be a faithful
representation of $H$ in $G = SL(V)$ where $dim(V) = r$. We shall use
the following notation when we wish to stress the fact that the quot
scheme as parametrising semistable principal quasibundles with
structure group $G$ rather than torsion--free sheaves of rank $r$:
\[
R_G := R^{\mu ss}
\]
Denote the universal quotient on $X \times R_G$ by $q_{_{R_G}} :
{\pi}_X^{*}{\co}_{X}(-m) \otimes W \lr {\ce}_G$.

Recall our notation of \S2 and the notion of quasibundles with respect
to the representation $\rho$. These are obtained from the affine
$X$--scheme $S_{\ce}$ together with a {\it generalised reduction}
datum $\sigma: X \lr S_{\ce}//H$. The pair $(\ce, \sigma)$ is a
principal quasibundle with structure group $H$. We remark that this
reduction datum $\sigma$ can be equivalently viewed as giving an
${\co}_X$--algebra morphism $\tau: (Sym^{*}({\ce} \otimes V))^{H} \lr
{\co}_X$. The map $\tau$ which come from a {\it genuine} generalised
reduction $\sigma$ is not simply the projection onto ${\co}_X$. It is
obtained by dualising the map $\sigma$, where ${\sigma}|_{U}$ is a
reduction of structure group of the principal bundle (associated to
${\ce}|_{U}$) to $H$ on a {\it large} open subset $U = U(\ce)$.

Let $R_H(\rho)$ be the scheme which parametrises pairs $(q,\tau)$
where $q \in R_G$ and $\tau$ is a reduction datum giving a {\it
quasibundle with structure group} $H$ and let $R_H(\rho)^o$ be the
open subset of pairs $(q,\tau)$ such that $\tau$ defines a {\it
principal bundle} with structure group $H$.

The existence of total families for principal quasibundles follows by
the general theory of Hilbert schemes and is shown in
(\cite{schmitt}); in Schmitt's language these will be parametrising
what he terms {\it honest singular bundles}. We have been somewhat
loose in defining the parameter space but we refer the reader to
Schmitt's paper (cf. \cite[Section 6.7]{schmitt} for details.

The scheme $R_H(\rho)$ is an $R_G$--scheme and the natural map
$f_{\rho}: R_H(\rho) \lr R_G$ is the one induced by $\rho$. Therefore,
since we have already fixed the Chern classes of the torsion free
sheaves in $R_G$, by the general theory of characteristic classes for
principal bundles, the characteristic classes of the principal bundles
in the open subscheme $R_H(\rho)^{0} \subset R_H(\rho)$ will also be
canonically fixed. This follows by the basic result of
Borel--Hirzebruch which relates the characteristic classes of
principal bundles with the Chern classes of associated vector
bundles. For a nice exposition and explicit results cf. \cite[Prop
3.2]{beauville}.

For simplicity we will denote this entire datum by ${\sf c}$ (we will
return to this in 6.2.1).

\brem It is immediate that the $\cg$-action on $R_G$ lifts to an
action on $R_H(\rho)$.  \erem

By the universal property of the scheme $R_H(\rho)$ it follows that
there exists a tautological family on $X \times R_H(\rho)$. Let
$\mathcal P$ denote this $R_H(\rho)$--flat family of $H$--quasibundles
on $X$ associated to the faithful representation $\rho$. 

Then by its definition, there exists an $R_H(\rho)$--flat family of
torsion free sheaves ${\mathcal P}(\rho) = {\cf}$ on $X \times
R_H(\rho)$. This is given precisely by the family $(id_X \times
f_{\rho})^*(\ce_G)$ obtained from the family of semistable torsion
free sheaves on $X \times R_G$. As remarked above, there is
$\cg$--action on $R_H(\rho)$ such that the family $\cf$ carries a
linearisation with respect to the $\cg$--action.

\bdefe\label{lambdadef} Let ${\mathcal P}(\rho) = {\cf}$ be the
induced family of torsion free sheaves on $X$ parametrised by
$R_H(\rho)$. Fixing the class $c \in K(X)$ as seen above, we have a
canonical choice of $u_1(c) \in K(X)$. We then have the {\it determinant
line bundles} $\mathcal L := {\lambda}_{\cf}(u_1(c))$ on $R_H(\rho)$
induced by the family $\cf$. \edefe

\noindent
{\bf Notation} {\it Since we have emphasized the role of the representation
$\rho$, we will henceforth denote the total space $R_H(\rho)$ simply
by $R_H$.}

Then we have the following:

\blem\label{huy} (\cite [Lemma 8.2.4]{hl})
\begin{enumerate}
\item[1.] If $s \in S$ is a point such that ${\mcp}_s|_{C}$ is
$\mu$--semistable. Then there exits an integer $n > 0$, and
$\cg$--invariant sections ${\overline \sigma}$ in
${\lambda}_F(u_1(c))^{n}$ such that ${\overline \sigma}(s) \neq 0$.

\item[2.] Let $s_1$ and $s_2$ be two points in $S$ such that either
${\mcp}_{s_1}|_{C}$ and ${\mcp}_{s_2}|_{C}$ are both semistable but
not $S$--equivalent (in the sense of Seshadri) or one of them is
semistable and the other is not. Then, in either case, there are
$\cg$--invariant sections in some tensor power of
${\lambda}_F(u_1(c))$ that separate $s_1$ and $s_2$.

\end{enumerate}
\elem

{\it Proof:} Since the proof follows the general line as given in the
text (\cite{hl}), we content ourselves by giving the main steps in the
argument and refer the reader to \cite{hl} for more details.

We need to relate the {\it ``determinant line bundle''} on $R_H$ to
its pull-back from the corresponding total family $Q_C$ of
principal $G$--bundles on a high degree curve $C \subset X$. 

Let $C$ be a general smooth curve in $X$ of high degree $a \gg 0$. Let
$i : C \hra X$. For the class $c \in K(X)$, its restriction, to
$c|_{C} := i^*(c)$ in $K(C)$, is completely determined by its rank $r$
since we have assumed that all our sheaves have trivial
determinant. Let $P'$ be the Hilbert polynomial determined on $C$ by
this restriction. Then for a large positive integer $m'$ we consider
the induced ``quot scheme'' of quotients $[q: W' \otimes {\co}_C(-m')
\lr E]$. Let us denote this quot scheme by $Q_C$. Here $W'$ is a
complex vector space of dimension $P'(m')$. Let $Q_C^{ss}$ be the
quotients which give semistable bundles on $C$ and let $q_C : W'
\otimes {\pi}_C^{*}({\co}_C(-m')) \lr {\ce}_C$ be the universal
quotient on $X \times Q_C$ with {\it determinant line bundle}
${\cl}_0'$ on $Q_C$ obtained from this family.  We then have the
following properties on {\it separating points} in $Q_C$ by sections
of the determinant line bundle ${\cl}_0'$:
\begin{enumerate}
\item All $SL(P'(m'))$--semistable points (in the GIT sense) in $Q_C$
are precisely the bundles $E$ which are semistable. Further, there
exists an integer $\nu$ and an $SL(P'(m'))$--invariant section
$\sigma$ of ${{\cl}_0'}^{\nu}$ such that ${\sigma}([q]) \neq 0$, where
$q \in Q_C^{ss}$ corresponds to the bundle $E$.

\item Two points $q_i$, $i = 1,2$ in $Q_C$ are separated by invariant
sections in some tensor power of ${\cl}_0'$ if and only if either both
are semistable but not $S$--equivalent or one of them is semistable and
the other is not.

\end{enumerate}

We now work with $R_H$ and the family ${\cf}$ on $X \times R_H$. Let
us denote the restriction of ${\cf}$ to $C \times R_H$ by ${\cf}_C$.

Now by increasing $m'$ if necessary, we can make sure that:
\begin{itemize}
\item The restrictions ${\cf}_s|_{C}$ are points of $Q_C$.

\item The sheaf $p_{*}{\cf_C}(m')$ is a locally free
${\co}_{R_H}$--module of rank $P'(m')$, where $p : C \times S \lr S$
is the projection.
\end{itemize}

Let $S_H$ be the associated projective frame bundle and let ${\eta}:
S_H \lr R_H$ be the natural map. It parametrises a quotient
\[
{\co}_{S_H} \otimes {\pi}_{C}^{*}{\co}_{C}(-m') \otimes W' \lr
{\eta}^{*}{\cf}_C \otimes {\co}_{\eta}(1).
\]
This induces a $SL(P'(m'))$--invariant morphism 
\[
{\Phi}_{\cf_{C}} : S_H \lr Q_C
\]
We also have the $\cg$--action on $\cf$ and via this action we get a
$\cg$--action on $S_H$ which commutes with the $SL(P'(m'))$--action
such that both $\eta: S_H \lr R_H$ and ${\Phi}_{\cf_{C}}$ are
equivariant for the ${\cg} \times SL(P'(m'))$--actions.

One then shows by a degree computation that:
\[
{{\Phi}_{\cf_{C}}^{*}({\cl}_0')}^{deg C} \simeq
{{\eta}^{*}{\cl}^{a^{2} deg(X)}} \eqno (\#)
\]
(The canonical {\it det} line bundle on $R_H$ depends only on the
choice of a polarisation on $X$ and the relation $\#$ is independent
of $C \in |a \Theta|$. cf. \cite{hl} pp 187-188)

Since $S_H \lr R_H$ is a principal $PGL(P'(m'))$--bundle, a section
invariant under $SL(P'(m'))$ descends to give a $\cg$--invariant
section. Thus, for any $\nu$ we thus get a linear map:
\[
s_{\cf}: {H^0(Q_C, {({\cl}_{0}')}^{\nu deg (C)})}^{SL(P'(m'))} \lr
{H^0(R_H, {\cl}^{\nu a^{2} deg(X)})}^{\cg}
\]
To conclude the proof of the Lemma we first observe that the
semistability of the quasibundle ${\mcp}_s$ is equivalent to the
$\mu$--semistability of the torsion--free sheaf ${\cf}_s$ for each $s
\in R_H$. Hence by using section of powers of the line bundle
${\cl}_0'$ on $Q_C$ and via the map $s_{\cf}$ we can indeed separate
the restrictions of quasibundles in $R_H$ to $C$ as well.

\begin{flushright} {{\it q.e.d}}\end{flushright}

\brem By the Mehta--Ramanathan theorem, if the degree of the curve $C$
is large enough, then the restriction of a $\mu$--semistable
torsion--free sheaf on $X$ to the curve $C$ is semistable. But the
same curve $C$ may not work for all the bundles in $R_H$. Hence we
need to work with the whole of $Q_C$. \erem

We have the following immediate corollary from the first part of Lemma
\ref{huy}:

\bcor\label{lepot} There exists an integer $\nu > 0$ such that the
line bundle ${\cl}^{\nu}$ on $R_H$ is generated by $\cg$--invariant
global sections i.e {\sf $\cl$ is $\cg$--{semi--ample}}. \ecor

Since $R_H$ is a quasi--projective scheme and since $\cl$ is
$\cg$--{\it semi--ample}, there exists a finite dimensional vector
space $J \subset J_\nu:= H^0(R_H,\cl^{\nu})^{\cg}$ that generates
$\cl^\nu$. Note that there is nothing canonical in the choice of $J$.

Let morphism $\phi_J: R_H \rightarrow \mathbb P(J)$ be the induced
$\cg$--invariant morphism defined by the sections in $J$. Because of
nonuniqueness of $J$, each choice of subspace of invariant sections
gives rise to a different map $\phi_{J'}$ to a different projective
space $\mathbb P(J')$.

\bdefe We denote the by $M_J$ the {\sf schematic image} $\phi_J(R_H)$
with the {\sf canonical reduced scheme structure}.  \edefe

\brem By the following result which may be titled $\cg$--properness, the
variety $M_J$ is {\it proper} and hence because of its
quasi--projectivity it is a projective variety. We note that we use
the term {\it variety} in a more general sense of an {\it reduced
  algebraic scheme of finite type} which need not be irreducible.  So
in what follows we will be working with the $\mathbb C$--valued points
of $M_J$.\erem

\bprop If $T$ is a separated scheme of finite type over $k$, and if
$\phi: R_H \lr T$ is an $\cg$ invariant morphism then
image of $\phi$ is proper over $k$.  \eprop 

{\it Proof}: This is an immediate consequence of Theorem \ref{ssred}
and \cite[Prop 8.2.6]{hl}. \begin{flushright} {{\it
      q.e.d}}\end{flushright}

\noindent
Let $J_\nu$ denote the vector space $H^0(R_H,\cl^{\nu})^\cg$, $\nu \in
\mathbb Z^{+};$ and Let $J \subset J_\nu$ be a finite dimensional
vector space which generates $\cl^\nu$.

For any $d \geq 1$, let $J^d$ be the image of the canonical
multiplication map $f_d : J \otimes,\cdots,\otimes J (d-fold)
\rightarrow J_{d\nu}$; in particular, $J^1 = J$.

Let $J'$ be any finite dimensional vector subspace of $J_{d\nu}$
containing $J^d$. Then clearly the line bundle $\cl^{d\nu}$ is also
globally generated by $\cg$--invariant sections coming from the subspace
$J'$ and this is so for any $d \ge 0$.  

So we have inclusions $J \hookrightarrow J^d \hookrightarrow J' $, and
hence a commutative diagram
\[
\xymatrix{
M_{J'} \ar[r]^{\pi_{J'/J}} & M_J \\
R_H \ar[u]_{\phi_{J'}} \ar[ru]_{\phi_J}
}
\]
Since $M_J$ and $M_J'$ are both projective, the map $\pi_{J'/J}$ is a
finite map. So if we fix a $J$ as above, we get an inverse system
(indexed by the $d \ge 1$) of projective varieties
$(M_{J'},\pi_{J'/J})$ and dominated by the finite type scheme $R_H$.

\[
\xymatrix{
& R_H \ar[dl] \ar[d]_{\phi_{J'}} \ar[dr]^{\phi_{J}} \\
\cdots \ar[r] & M_{J'} \ar[r]_{\pi_{J'/J}} & M_{J}  
}
\]

Hence the inverse limit of the system $(M_{J'},\phi_{*})$ is in fact one of
the $M_{J'}$'s where $J'$ is a finite dimensional subspace of
$H^0(R_H,\cl^{n})^{\cg}$ which generates $\cl^n$.  

\bdefe We denote this {\it inverse limit variety} by $M_H(\rho)$ and
let $\pi:R_H \rightarrow M_H(\rho)$ be the canonical morphism induced
by the invariant sections coming from the subspace $J_0$
associated to the inverse limit.\edefe

We summarise the above discussion in following theorem:

\bth\label{lepo} Then the reduced projective scheme $M_H(\rho)$
parametrises equivalence classes of $\mu$--semistable quasibundles
with structure group $H$.  There is a natural morphism $q : R_H \lr
M_H(\rho)$. Furthermore, the representation $\rho: H \lr G = SL(V)$
induces a natural morphism ${\overline f_{\rho}} : M_H(\rho) \lr M_G$,
where $M_G := M_{SL(V)}$ is the moduli space of the $\mu$--semistable
torsion free sheaves with trivial determinant and fixed $c_2$. \eeth

{\it Proof:} The existence of the morphism ${\overline f_{\rho}} :
M_H(\rho) \lr M_G$ follows by the {\it naturality} of the moduli space
$M_G$ by virtue of the existence of the family of semistable torsion
free sheaves $\cf$ on $R_H(\rho)$.

\begin{flushright} {\it q.e.d} \end{flushright}

\brem The strategy is somewhat similar to that of \cite{basa} but
unlike \cite{basa}, there are no GIT quotients involved here. \erem

\brem Note that this is not a categorical quotient since $\cl$ is not
ample and is only {\it semi-ample} (Cor \ref{lepot}), i.e some power
of $\cl$ is generated by sections.  \erem

\subsection{Towards the description of the moduli space}

In order to get a better understanding of the geometry of the moduli
space $M_H(\rho)$ we need the following lemmas:

\blem\label{polyeq} Let $P_1$ and $P_2$ be two polystable principal
$H$--bundles on $X$. Let $a \gg 0$ and $C \in |a \Theta|$ be a general
smooth curve. Then $P_1 \simeq P_2$ if and only if $P_1|_{C} \simeq
P_2|_{C}$.  \elem

{\it Proof:} Consider the principal $H \times H$--bundle $P_1 \times_X
P_2$. Let us denote this bundle by $E$.

By (\cite{serre}, p 19), to give an isomorphism between $P_1$ and $P_2$
is equivalent to giving a {\it reduction of structure group} of $E$ to
the diagonal embedding $\Delta \subset (H \times H)$.

By the Chevalley semi--invariant theorem, we can embed $(H \times
H)/{\Delta} \subset W$ in a $H \times H$--module $W$ as a closed
orbit.

Thus a section of $E(\frac{H \times H}{\Delta})$ is a section of the
vector bundle $E(W)$ which lies in $E(\frac{H \times H}{\Delta})$.

\noindent
{\underline{\it Claim}}: The principal bundle $E$ is also polystable
of degree $0$. Let us assume this claim and complete the proof.

By the usual Enriques--Severi lemma for a general high degree curve
$C$, the restriction map:
\[
H^0(X, E(W)) \lr H^0(C, E(W)|_{C}) \eqno(\&)
\]
is {\it surjective}. Further, by assumption, for high degree curves we
know that $P_1|_{C} \simeq P_2|_{C}$. Hence, by our discussion above
we have a reduction of structure group of $E|_{C}$ to $\Delta$. i.e a
section:
\[
\sigma_{C}: C \lr E(\frac{H \times H}{\Delta})|_{C}~~\hra E(W)|_{C}
\]
By $(\&)$ above, there exists a lift of $\sigma_C$ to a section
$\sigma: X \lr E(W)$. We need to show that the image of $\sigma$ lies
in $E(\frac{H \times H}{\Delta})$.

One knows that $E$ is polystable of degree $0$, by the {\it Claim}
above.  Thus we are in the setting of Prop \ref{lcs}. Thus we
can get a big open subset $U \subset X$ such that the image of the
section $\sigma(U) \subset E(\frac{H \times H}{\Delta})$.

Now since $\frac{H \times H}{\Delta} \subset W$ is a closed embedding, it
follows that the entire image $\sigma(X) \subset E(\frac{H \times
H}{\Delta})$. This gives the required reduction of structure group to
$\Delta$. To complete the proof we need to prove the claim.

To see this we again use the easier half of the Mehta--Ramanathan
restriction theorem (i.e Lemma \ref{rest}). Thus to show that $E$ is
polystable of degree $0$ we need to show that for a general high
degree curve $C$, the restriction $E|_{C}$ is polystable of degree
zero.

Since $P_1|_{C}$ and $P_2|_{C}$ are polystable of degree $0$ (by the
restriction theorem again), it follows by Ramanathan's theorem
(\cite{rama}) (the Narasimhan--Seshadri theorem for principal bundles
on curves) that there exists representations $\rho_i : {\pi}_1(C) \lr
K$ $i = 1,2$, $K$ a maximal compact subgroup of $H$, such that $P_i|_{C}$
is the $H$--bundle associated to $\rho_i$, for $i = 1,2$.

It is easy to see that the bundle $E|_{C}$ is the bundle associated to
the representation $\rho_1 \times \rho_2$ (\cite{rama} p 146) and the
polystability follows by Ramanathan's theorem again. This completes
the proof of the Lemma.
\begin{flushright} {\it q.e.d} \end{flushright} 

%\brem We prove a more general fact in the appendix.\erem

\subsubsection{Associated graded of a semistable quasibundle}
Let ${\cp} \lr X$ be an $H$--quasibundle. Let $\rho : H \hra GL(V)$ be
the accompanying faithful representation and ${\ce} = {\cp}(\rho)$ the
associated torsion--free sheaf. The sheaf $\ce$ possesses a
Jordan--Holder filtration ${J}^{\bullet}$ by saturated subsheaves and
we may take the associated graded sheaf
$gr_{_{J^{\bullet}}}(\ce)$. For two different filtrations
$J_i^{\bullet}$ we have an isomorphism of vector bundles
$(gr_{_{J_1^{\bullet}}}(\ce))^{**} \simeq
(gr_{_{J_2^{\bullet}}}(\ce))^{**}$ By an abuse of notation we write
$gr^{\mu}(\ce)$ to denote {\it an} associated graded torsion--free
sheaf.

Let ${\ba} = {\ba}^{1}_{\bc}$. Recall that there exists a family
$\{{\ce}_t\}_{t \in \ba}$ of torsion--free sheaves such that ${\ce}_t
\simeq {\ce}$ for $t \neq 0$ and ${\ce}_0 \simeq gr^{\mu}(\ce)$.

By the {\it Semistable reduction Theorem \ref{ssred} (and its proof)}
we have an $H_A$--quasibundle (where $Spec A = {\ba}_1 \subset {\ba}$
an open subset which contains $0$). Denote this $H_A$--quasibundle by
$\{{\cp}_t\}_{t \in {\ba}_1}$. Further, this family has the property
that ${\cp}_t \simeq {\cp}$ for $t \neq 0$ and there is a big open set $U$ 
and a rational $H'$--bundle ${\cp}_0'$ on $U$, such that ${\cp}_0|_{U}$ 
is obtained from ${\cp}_0'$ by an extension of structure group
via $\psi: H' \lr H$ (recall that $H'$ could be a non-reductive group; cf.
Prop \ref{serre} and Lemma \ref{ssnonred2}).

Also from the proof (Theorem \ref{flatext}) one knows that
${\cp}_0'(\rho) \simeq {\ce}_0|_{U}$ which is {\it polystable}.

{\underline{\it Claim}}: The quasibundle ${\cp}_0$ is also polystable.

By Lemma \ref{rest} it is enough to show that ${\cp}_0|_{C}$ is
polystable, where $C$ is a high degree curve contained in the big open
subset $U$, where ${\cp}_0$ is a principal bundle. Further,
${\cp}_0'|_{C}$ is a bundle which is {\it flat} i.e it comes from a
representation $\chi_1 : {\pi}_1(C) \lr H'$ (see Prop \ref{lc}
). Moreover, ${\cp}_0'(\rho)|_{C} \simeq {\ce}_0|_{C}$ which is {\it
polystable} by the Mehta--Ramanathan theorem, since ${\ce}_0|_{U}$ is
polystable. In other words, under our definition of polystability of
principal bundles with non--reductive structure groups (i.e Def
\ref{ssnonred1}), the bundle ${\cp}_0'|_{C}$ is {\it polystable}.

Again ${\cp}_0|_{C}$ comes from ${\cp}_0'|_{C}$ by extension of
structure group $f: H' \lr H$. Therefore, since ${\cp}_0'|_{C}$ is
polystable by Prop \ref{ssnonred}, it follows that ${\cp}_0|_{C}$ is
{\it polystable}. This proves the claim.

\bdefe\label{assgr} Let $\cp$ be a semistable $H$--quasibundle. If
there exists a family $\{{\cp}_t\}_{t \in Spec A}$, with $A$ a
complete discrete valuation ring, such that ${\cp}_t \simeq {\cp}$,
for $t \neq 0$ and ${\cp}_0$ polystable. Then we call ${\cp}_0$ {\bf
an} associated graded quasibundle of ${\cp}$. \edefe

This is uniquely defined upto {\it double duals} in the following
sense:

\blem Let ${\cp}$ be an $H$--quasibundle. Then if ${\cp}_{s_0}$ and
${\cp}_{t_0}$ are two choices of polystable limits (as above), then
${\cp}_{s_0}^{**} \simeq {\cp}_{t_0}^{**}$.  \elem

{\it Proof:} Let ${\cp}_S$ and ${\cp}_T$ be two families of
quasibundles such that at the closed points $s_0 \in S$ and $t_0 \in T$
the quasibundles ${\cp}_{s_0}$ and ${\cp}_{t_0}$ are polystable and
the generic fibres in either family is isomorphic to $\cp$.

Consider the open subset $U \subset X$ where both ${\cp}_{s_0}$ and
${\cp}_{t_0}$ are locally free. Then since $S$ and $T$ are spectra of
discrete valuation rings, the families ${\cp}_S$ and ${\cp}_T$ are
locally free on $U \times S$ and $U \times T$ respectively. Using the
Mehta--Ramanathan theorem choose a general high degree curve $C
\subset U$ so that the restrictions ${\cp}_{s_0}|_{C}$ and
${\cp}_{t_0}|_{C}$ are polystable. Since $C \subset U$ it follows that
both ${\cp}_{s_0}|_{C}$ and ${\cp}_{t_0}|_{C}$ are polystable limits
of ${\cp}|_{C}$ (which is also semistable by openness of
semistability). But the associated graded of a semistable principal
bundle on a smooth projective curve is uniquely defined (cf.
\cite{r1}). Hence,
\[
{\cp}_{s_0}|_{C} \simeq {\cp}_{t_0}|_{C}.
\]
Now we apply Lemma \ref{polyeq} to the double dual principal bundles
on $X$ and we are done.

\brem If ${\cp}(\rho) = {\ce}$, then for each choice of associated
$H$--quasibundle $gr^{\mu}({\cp}) = {\cp}_0$ defined above, we see
that ${\cp}_0^{**}(\rho) \simeq {\ce}_0^{**}$. \erem

\section{The geometry of the moduli of $H$--bundles}

{\em From this section onwards $H$ is a simple algebraic group}. In
this section we study the points of the moduli space ${\ol {M^{0}_H}}$
intrinsically as well as in relationship with its image points in the
moduli space $M_G$.

\subsection{Cycles associated to quasibundles in ${\ol {M^{0}_H}}$}

Let $\cp$ be a semistable $H$--quasibundle in ${\ol
{M^{0}_H}}$. Therefore there exists a family $\{{\cp}_t\}_{t \in Spec
A}$, with $A$ a complete discrete valuation ring, such that ${\cp}_t$
is a semistable principal $H$--bundle, for $t \neq 0$ and ${\cp}_0
\simeq {\cp}$. 

Our aim is to associate a cycle $Z_{\cp} \in S^l(X)$ of degree $l$ to
the quasibundle $\cp$ in an intrinsic manner. By this we mean that the
pair $({\cp}^{**}, Z_{\cp})$ in the compactified moduli space is
independent of the quasibundle (and hence the representation $\rho$ as
well). Towards this we work with the chosen representation $\rho$ used
in defining the quasibundle $\cp$.

\subsubsection{The Dynkin index of the representation $\rho$}

We recall the notion of Dynkin index (\cite{dynkin}, \cite{snr},
\cite{atiyah}) and some basic results from these sources which will
play a key role in what follows.

\bdefe\label{dynkin} Let $\theta$ be a the highest root in ${\goth H}$
the Lie algebra of $H$ and ${\goth sl(\theta)}$ the $3$--dimensional
sub--algebra of ${\goth H}$ associated to it. Decompose the
$H$--module $V$ into ${\goth sl(\theta)}$--modules as ${\oplus V_i} =
V$. Let $dim(V_i) = m_i$. Then we can define the Dynkin index
$m_{\rho}$ of the $H$--module $V$ as follows:
\[
m_{\rho} = \sum_{i} \binom{m_i+1}{3}
\]
\edefe

\brem The Dynkin index of a simple subgroup $H$ of a simple Lie group
$G$ is usually defined as the ratio of the invariant inner product on
$Lie(H)$ to the invariant inner product of $Lie(G)$ where the inner
products are normalised to make the length of the highest root $2$
(\cite[page 455]{atiyah}). \erem

\beg
\begin{enumerate}

\item Let ${\rho}_m : SL(2) \lr SL(V)$ be the standard irreducible
representation $V = S^{m}(W)$ where $dim(W) = 2$. The $m_{\rho_{m}} = 
\binom{m+2}{3}$.

\item Let ${\rho} : SL(2) \lr SL(V)$ be the inclusion obtained by
identifying $SL(2)$ with the $3$--dimensional group given by the
highest root $\theta$ in ${\goth sl(V)}$. Then $m_{\rho} = 1$.

\item If $\rho$ is the adjoint representation of $H$ then the index is
$2 h$, where $h$ is the dual Coxeter number of $H$.
\end{enumerate}
\eeg

The following is the list of index $1$ subgroups given by Dynkin
(\cite{dynkin}). We reproduce the list from \cite{atiyah}:

\begin{center}
\tiny
TABLE 1

$S(U(n) \times U(1)) \subset SU(n+1)$

$Spin(n) \subset Spin(n+1)$

$Sp(n) \times Sp(1) \subset Sp(n+1)$

$SU(3) \subset G_2$

$Spin(9) \subset F_4$

$F_4 \subset E_6$

$E_6 \times U(1) \subset E_7$

$E_7 \times SU(2) \subset E_8$

\end{center}

\normalsize

We now quote from \cite{snr} (see also \cite[page 455]{atiyah});

\bprop\label{pythree} Let ${\rho}_{*} : {\pi}_3(H) \lr
{\pi}_3(SL(V))$. Then it is given by {\em ``multiplication''} by
$m_{\rho}$. \eprop

\brem\label{zeroscheme} This in particular proves that any index $1$
inclusion $SU(2) \subset SU(n)$ preserves the second Chern class of a
principal $SU(2)$--bundle when we extend structure group to $SU(n)$
since ${\pi}_3(G) \simeq {\bz}$ classifies principal $G$--bundles on a
real $4$--sphere $S^4$ and the second Chern class of the bundle
classifies it. This also implies that any $SU(n)$--bundle $E$ on $S^4$
is obtained as the extension of structure group from $SU(2)$. In
particular, the the rank $n$ vector bundle underlying $E$
(topologically) splits as a direct sum of a rank $2$ bundle and direct
sum of $n - 2$ trivial line bundles. Thus if the vector bundle
underlying $E$ has a nontrivial section $s$ then the section comes
from a section of the rank $2$ subbundle. Hence the {\it
``zero--scheme''} $Z(s)$ is represented by $c_2(E)$. \erem

For a point $x \in X$, let $B_{\epsilon}(x)$ be an analytic ball of
radius $\epsilon$ around $x$. Let $E$ be a vector bundle on $X$ of
rank $r$ and let $E_{|_{B_{\epsilon}(x)}}$ be its restriction to
$B_{\epsilon}(x)$. Suppose further that it is trivial on the boundary
$\partial{B_{\epsilon}(x)}$. Let us denote by $E_{|_{B_{\epsilon}(x)}}
/E_{|_{\partial{B_{\epsilon}(x)}}}$ the complex vector bundle on the
real four sphere $S^4 \simeq B_{\epsilon}(x)
/{\partial{B_{\epsilon}(x)}}$ obtained by identifying
$E_{|{\partial{B_{\epsilon}(x)}}}$ to ${\bc}^{\oplus r}$ using the
trivialisation.

 We then have the following result from
\cite{li} and \cite{morgan} extended for rank $r$ bundles. Li and
Morgan show it for rank $2$ case but the proofs generalise to the
higher rank case.

\bprop\label{sphere} Let ${\ce}_{\ba}$ be a family of (analytic)
torsion--free sheaves of rank $r$ over $B \times {\ba}$ flat over
$\ba$, where $B$ is a closed ball in the surface $X$ and $\ba$ is the
affine line. Further let ${\ce}_{\ba}$ be locally free on $B \times
{\ba} - 0$ and ${\ce}_0$ be torsion--free with singularity at the
origin $0 \in B$. Let $B_{\epsilon,u} \subset B \times \{u\}$ be the
2--dimensional $\epsilon$--ball centered at $(0,u)$ and let ${\ce}_u$
be the restriction of $\ce$ to $B_{\epsilon,u}$. Then the smooth
trivialisation of ${\ce}_0{_{|_{\partial B_{\epsilon,u}}}}$ which is
induced from a trivialisation of
$({\ce}_0{_{|_{B_{\epsilon,u}}}})^{**}$ induces a family of
trivialisations $\beta_u: {\bc}^{\oplus r} \times {\partial
B_{\epsilon,u}} \simeq {\ce}_u{_{|_{\partial B_{\epsilon,u}}}}$ which
depends smoothly on $u$ with $u \in {\ba}$ being
$\epsilon$--small. Furthermore, we also have on $S^4 \simeq
B_{\epsilon,u} /{\partial{B_{\epsilon,u}}}$ the following:
\[
c_2({\ce}_u{_{|_{B_{\epsilon,u}}}}
/{\ce}_u{_{|_{\partial{B_{\epsilon,u}}}}}) = length ({\ce}_0^{**} /
{\ce}_0)_{0}
\]
\eprop

{\it Proof}: The only new ingredients needed to generalise Li's
arguments are:
\begin{itemize}
\item Extension of \cite[Prop 6.4]{li} to higher ranks which has been
done by Gieseker and Li \cite{giesli}.

\item When the length $l = 1$ then the length $l$ can be realised as
the second Chern class of the rank $r$ bundle on $S^4$. This is done
in Li in the rank $2$ case by getting a section whose zero-scheme has
length $1$ and is represented by the $c_2$. Since we are on $S^4$ this
argument goes through for higher rank case as well as we have seen in
Remark \ref{zeroscheme}.

\end{itemize}
\begin{flushright} {\it q.e.d} \end{flushright}

Let ${\cp}_{t_0} \in {\ol {M^{0}_H}}$ be a polystable $H$--quasibundle
and ${\cp}_T$ be a family of semistable principal $H$--bundles with
${\cp}_{t_0}$ as limit. Let ${\ce}_T = {\cp}_T(\rho)$. Let
$Z({\ce}_{t_0})$ be the cycle associated to the torsion--free sheaf
${\ce}_{t_0}$. Recall that $Z({\ce}_{t_0}) \in S^l(X)$ where $l =
c_2({\ce}_{t_0}^{**}) - c_2({\ce}_{t_0})$. Further the cycle is given
by:
\[
Z({\ce}_{t_0}) := \sum_{x \in X} length({\ce}_{t_0}^{**}/{\ce}_{t_0})_{x} \cdot x
\]

\bth\label{dynkimulti} Let the notations be as above. Let $x \in
Supp(Z({\ce}_{t_0}))$. Then the number
$length({\ce}_{t_0}^{**}/{\ce}_{t_0})_{x}$ is a multiple of the Dynkin
index $m_{\rho}$. In particular, the total degree $l$ of the cycle
$Z({\ce}_{t_0})$ is also a multiple of $m_{\rho}$. \eeth

{\it Proof}: Fix a point $x \in Supp(Z({\ce}_{t_0}))$ and choose an
ball $B$ around it and identify $x$ with $0 \in B$. Let $u \in {\ba}$
is the point in a small disk in $\ba$ corresponding to $t_0 \in T$.

By restricting ${\ce}_T$ to $B \times {\ba}$, in the notation of Prop
\ref{sphere}, we have a family ${\ce}_{\ba}$ satisfying the properties
given there. Furthermore, we also have the extra datum that the
bundles ${\ce}_t$ for $t \neq t_0$ have reduction of structure group
to principal $H$--bundles ${\cp}_t$ $\forall t \neq t_0$.

Following the procedure in Prop \ref{sphere} we get a vector bundle
${\ce}_u{_{|_{B_{\epsilon,u}}}}
/{\ce}_u{_{|_{\partial{B_{\epsilon,u}}}}}$ such that
\[
c_2({\ce}_u{_{|_{B_{\epsilon,u}}}}
/{\ce}_u{_{|_{\partial{B_{\epsilon,u}}}}}) = length ({\ce}_0^{**} /
{\ce}_0)_{0}
\]

We now observe easily that the vector bundle
${\ce}_u{_{|_{B_{\epsilon,u}}}}
/{\ce}_u{_{|_{\partial{B_{\epsilon,u}}}}}$ on $S^4$ also has a
reduction of structure group to $H$. Thus we have a principal
$H$--bundle $P_{\epsilon}$ on $S^4$ whose extension of structure group
by $\rho: H \hra SL(V)$ is ${\ce}_u{_{|_{B_{\epsilon,u}}}}
/{\ce}_u{_{|_{\partial{B_{\epsilon,u}}}}}$.

By Prop \ref{pythree} it follows that the second Chern class of
${\ce}_u{_{|_{B_{\epsilon,u}}}}
/{\ce}_u{_{|_{\partial{B_{\epsilon,u}}}}}$ is a multiple of the Dynkin
index $m_{\rho}$. Hence $length ({\ce}_0^{**} / {\ce}_0)_{0}$ is a
multiple of $m_{\rho}$. This proves the theorem. \begin{flushright}
{\it q.e.d} \end{flushright}

\bdefe\label{cycle} Let ${\cp}_{t_0} \in {\ol {M^{0}_H}}$. The we define the cycle
$Z({\cp}_{t_0})$ as:
\[
Z({\cp}_{t_0}) := \sum_{x \in X}
\frac{length({\ce}_{t_0}^{**}/{\ce}_{t_0})_{x}}{m_{\rho}} \cdot x
\]
\edefe

\subsection{\it Points of the moduli}

Let ${\cp}$ be an $H$--quasibundle and let ${\cp}(\rho) =
{\ce}$. Further, let ${\ce}_0$ be an associated graded sheaf of
${\ce}$ and let us denote the double dual ${\ce_0}^{**}$ by $\be$.

Let $Q({\be},{\it l})$ be the quot scheme of torsion quotients of $\be$
of length ${\it l}$. By the choice of our $\be$, the torsion--free
sheaves $F \in Q({\be},{\it l})$ have all {\it reduction of structure
groups} to quasibundles. 

{\it Notation}: Let the induced total family quasibundles coming from
$Q({\be},{\it l})$ be denoted by $Q_{H}({\be},{\it l})$. We thus have
the following diagram:

\begin{center}
\begin{picture}(400,325)(0,20)
\put(-8,290){\ctext{$Q_H({\be},{\it l})$}}
\put(75,290){\vector(1,0){200}}
\put(349,290){\ctext{$Q({\be},{\it l})$}}
\put(325,250){\vector(0,-1){200}}
\put(25,250){\vector(0,-2){200}}

\put(25,25){\ctext{$S^{\it l_{0}}(X)$}}
\put(95,25){\vector(1,0){175}}

\put(75,265){\vector(1,-1){200}}
\put(325,25){\ctext{$S^{\it l}(X)$}}
\put(125,150){\ctext{$h$}}
\put(350,150){\ctext{$q$}}
\put(0,150){\ctext{$h_0$}}
%\put(250,200){\ctext{$\circlearrowright$}}
\put(175,320){\ctext{$f_{\rho}$}}
\put(175,50){\ctext{$''m_{\rho}''$}}

\end{picture}
\end{center}
We note the following:
\begin{enumerate}
\item The degree $l_0$ is therefore $l_0 = {l \over m_{\rho}}$.
\item The vertical map $h_0: Q_H({\be}, {\it l}) \lr S^{\it
l_{0}}(X)$ is the one which associates to each ${\cp}$ the cycle
$Z({\cp})$ defined above.

\item The map $''m_{\rho}'': S^{\it l_{0}}(X) \lr S^{\it l}(X)$ is the map
induced by multiplication by the Dynkin index $m_{\rho}$ which maps
the cycle $Z({\cp})$ to $m_{\rho} \cdot Z({\cp}) =
Z({{\ce}_0})$. Observe that this map is an {\it injection}.
\end{enumerate}

\brem\label{action} {\it(Action of $Aut({\be})$ on the fibre of
$f_{\rho}$)}

Let $\alpha \in Aut({\be})$. Let ${\ce}$ be a ${\mu}$--polystable
torsion--free sheaf with ${\be} = gr^{\mu}(\ce)^{**}$. Let $\cp$ be a
quasibundle obtained from a reduction datum $\sigma: X \lr
S_{\ce}//H$.

We first observe that since $\be$ is a locally free sheaf on $X$, by
Hartogs theorem we see that the natural restriction map $res_U: Aut({\be})
\lr Aut({\be}_U)$ is an isomorphism. (this holds when $\be$ is reflexive, cf.
for example Cor 1.11.1 \cite{maruyama}).

By restricting $\alpha$ to $U$ a big open set, we have an action of
$\alpha$ on $\sigma_U$ (restriction of $\sigma$ to $U$). Let $\alpha
\cdot \sigma_U = \sigma'_U$. Then by definition $\sigma'_U: U \lr
S_{\ce}//H$ which extends {\it uniquely} to a new datum $\sigma': X
\lr S_{\ce}//H$. We define:
\[
\alpha \cdot \sigma = \sigma'
\]
It is clear that if ${\cp}'$ is the quasibundle obtained from
$\sigma'$, then ${\cp}'(\rho) = {\ce}$ and ${\cp}^{**} \simeq {\cp}'^{**}$.

\erem

We then have the following:

\bprop\label{connected} Let $\cp$ and $\cq$ be two $\mu$--polystable
$H$--quasibundles in $Q_H({\be},{\it l})$ lying in the same fibre of
$f_{\rho}$ . If ${\cp}^{**} \simeq {\cq}^{**}$ then $\cp$ and $\cq$
lie in the same connected component of the fibre of the map $h:
Q_H({\be},{\it l}) \lr S^{\it l}(X)$ (or equivalently of the map $h_0:
Q_H({\be},{\it l}) \lr S^{\it l_{0}}(X)$).  \eprop

{\it Proof:} We first show that if $\cp$ and $\cq$ are two {\it
polystable} quasibundles in a fibre of the map $f_{\rho}: Q_H({\be},{\it
l}) \lr Q({\be},{\it l})$ and such that ${\cp}^{**} \simeq {\cq}^{**}$
then $\cp$ and $\cq$ lie in the same orbit of the action $Aut({\be})$.

Note firstly that the polystability of $\cp$ and $\cq$ as well as the
property that they lie in the same fibre of $f_{\rho}$ forces that the
{\it associated torsion--free sheaves} in $Q({\be},{\it l})$ are
$\mu$--polystable Hence, both $\cp$ and $\cq$ are {\it reductions} of
structure groups of the ``same'' polystable torsion--free sheaf
${\ce}_0$. One also sees that ${\ce}_0^{**} \simeq {\be}$. In other
words, ${\cp}, {\cq} \in f_{\rho}^{-1}({\ce}_0)$.

Consider the restriction to a big open subset $U = U({\ce}_0)$. By
definition, both ${\cp}_U$ and ${\cq}_U$ are reductions of structure
group of ${\ce}_U \simeq {\be}_U$ and since ${\cp}^{**} \simeq
{\cq}^{**}$, these are isomorphic reductions of structure group of the
principal bundle ${\be}_U$. Therefore they are in the same orbit of
$Aut({\be}_U)$. By the definition of the action in Remark \ref{action},
this implies that $\cp$ and $\cq$ lie in the same orbit of $Aut({\be})$.

Now since $\be$ is a polystable vector bundle on $X$ (we confuse the
principal bundle with the vector bundle associated to it) by an
argument similar to the one in Lemma \ref{polyeq}, it follows that for a
general curve $C \in X$, $Aut({\be}) \simeq Aut({\be}|_{C})$. Also by the
Mehta-Ramanathan restriction theorem, ${\be}|_{C}$ is also a polystable
vector bundle on $C$ and hence its automorphism group is a product of
$GL(n)$'s. Hence, $Aut({\be})$ is connected. This implies that the orbit
of $Aut({\be})$ is connected and hence both $\cp$ and $\cq$ lie in the
same connected component of the fibre of $h$. 
\begin{flushright} {\it q.e.d} \end{flushright}

\bcor\label{keycor} Let $\cp$ and $\cq$ be two $\mu$--semistable
$H$--quasibundles in $Q_H({\be},{\it l})$ lying in the same fibre of
$f_{\rho}$. If ${\cp}_0^{**} \simeq {\cq}_0^{**}$ then $\cp$ and $\cq$
lie in the same connected component of the fibre of the map $h:
Q_H({\be},{\it l}) \lr S^{\it l}(X)$.\ecor

{\it Proof:} Observe that for any quasibundle $\cp$, an associated
graded quasibundle is connected to $\cp$ by a path (cf. Def
\ref{assgr}). Hence the Corollary follows from the Prop
\ref{connected}.

\bprop\label{keycor1} Let $\cp$ and $\cq$ be two $\mu$--polystable
$H$--quasibundles in $Q_H({\be},{\it l})$ lying in the same fibre of
$h$. If further, ${\cp}^{**} \simeq {\cq}^{**}$, then $\cp$ and $\cq$
lie in the same connected component of the fibre of the map $h:
Q_H({\be},{\it l}) \lr S^{\it l}(X)$ (or equivalently of the map
$h_0: Q_H({\be},{\it l}) \lr S^{\it l_{0}}(X)$).\eprop

{\it Proof:} By the result of of Baranovsky (\cite{bara}) and
Ellingsrud and Lehn (\cite{el}) (and in rank 2, Li (\cite[Prop
6.5]{li})), one knows that the fibre of the map $q$, which we denote
by $T$, is {\it connected}, in fact {\it irreducible}. In particular,
every torsion--free sheaf $\cf$ in the fibre of $q$ is a sheaf such
that ${\cf}^{**} \simeq {\be}$ and since $\be$ is polystable, it
follows that so is $\cf$. In other words, every closed point of the
fibre $T$ is represented by a polystable torsion--free sheaf.

Let us fix the polystable quasibundle $\cp$ in the fibre of $h$ and
let $\cq$ be another polystable quasibundle such that ${\cq}^{**}
\simeq {\cp}^{**}$. The polystable torsion--free sheaves {\it
associated} to $\cp$ and $\cq$, namely $\ce$ and $\cf$ are in $T$. If
$\ce$ and $\cf$ coincide then the result follows from Prop
\ref{connected}. Let $f_{\rho}(\cp) = t_0$ and $f_{\rho}(\cq) = t_1$.
Since $T$ is irreducible, without loss of generality we may take $T$
to be an irreducible curve joining $t_0$ and $t_1$. Let ${\ol T}$ be
its normalisation. Hence, ${\ol T}$ is a smooth irreducible curve.

{\it We {\sf claim} that $\cq$ lies in the same connected component of
the fibre of $h$ which contains $\cp$}.

Consider the family of torsion--free sheaves ${\ce}_T$ (on $X \times
T$) such that ${\ce}_{t_0} \simeq {\ce}$ and ${\ce}_{t_1} \simeq
{\cf}$, which exists by the definition of the Quot scheme of torsion
quotients. We shall work with the pull--back of ${\ce}_T$ to $X \times
{\ol T}$ and {\it denote it by ${\ce}_{\ol T}$}. Thus we may
take $t_0,t_1 \in {\ol T}$.

Since all the sheaves ${\ce}_t$ lie over the same fibre $T$ it follows
that they are all polystable sheaves and have the same double dual
$\be$. Let $\cu$ be the maximal open subset where ${\ce}_{\ol T}$ is locally
free. Then we have an isomorphism of families on $\cu$:
\[
{\ce}_{\ol T}|_{\cu} \simeq {\be}_{\ol T}|_{\cu}
\]
where ${\be}_{\ol T} = {\be} \times {\ol T}$.

Since the quasibundle $\cp$ is a ``reduction of structure group'' of
${\ce}_{t_0}$, it implies that there is a section ${\sigma}: X \lr
(S_{{\ce}_{t_0}})//H$.

Restricting this section to ${\cu}_{t_0} = {\cu} \cap (X \times
{t_0})$ we get a section ${\sigma}_{{\cu}_{t_0}} : {\cu}_{t_0} \lr
(S_{{\ce}_{t_0}})//H$. Note that ${\cu}_{t_0}$ is also a {\it big}
open subset of $X$ (cf. \cite[2.6 ,page 1187]{schmitt}).

Observe that the double dual principal bundle ${\cp}^{**}$ of $\cp$
is also obtained by a reduction of structure group of the principal
bundle associated to $\be$. 

Let us denote this reduction by $\tau: X \lr (S_{\be})/H$. On the big
open set ${\cu}_{t_0}$ this reduction ${\tau}_{{\cu}_{t_0}}$ and the
reduction $\sigma_{{\cu}_{t_0}}$ giving $\cp$ are therefore mapped to
each other by an automorphism $\phi \in Aut({\be}|_{{\cu}_{t_0}}) =
Aut(\be)$.

Since ${\be}_{\ol T} = {\be} \times {\ol T}$, the reduction section
$\tau$ trivially extends to a section $\tau_{\ol T} : X \times {\ol T}
\lr (S_{{\be}_{\ol T}})/H$. Restricting $\tau_{\ol T}$ to $\cu$ we get
a section
\[
{\tau}_{\cu} : {\cu} \lr (S_{{\be}_{\ol T}|_{\cu}})//H \simeq
(S_{{\ce}_{\ol T}})//H.
\]
Since $\cu$ is {\it big}, and since $X \times {\ol T}$ is {\it smooth
(and hence normal)}, this section extends to give a new section
${{\ol{\tau}}}: X \times {\ol T} \lr (S_{{{\ce}_{\ol T}}|_{X \times
{\ol T}}})//H$, i.e a ``reduction datum'' for the family ${\ce}_{\ol
T}$ (see Remark \ref{bigopen}).

This gives a family ${{\ol {\cp}}}_{\ol T}$ of quasibundles on $X
\times {\ol T}$ such that ${{\ol {\cp}}}_{t_0}$ and ${\cp}$ are in the
same orbit of $Aut(\be)$ (by the element $\phi$). Thus ${{\ol
{\cp}}}_{t_0}$ and $\cp$ lie in the same connected component of the
fibre of $h$.
 
By the universal property of $Q_H({\be},{\it l})$, we get a morphism
induced by this family of quasibundles, namely
\[
\psi: {\ol T} \lr Q_H({\be},{\it l})
\]
($\psi$ depends on the choice of the reduction $\tau$).

Further, by the definition of the reduction datum on $\cu$ all the
quasibundles in the family ${{\ol {\cp}}}_{\ol T}$ (which are polystable)
have {\it isomorphic double duals} ${\cp}^{**}$.

In particular, they lie in the same connected component of the fibre
of $h$. Moreover, if $\cq$ corresponds to the point $t_1 \in {\ol T}$
then the quasibundles ${{\ol {\cp}}}_{t_1}$ and $\cq$ lie in the same
fibre $f_{\rho}^{-1}({\ce}_{t_1})$. By Prop \ref{connected} they lie
in the same orbit of $Aut({\be})$ since ${\cq}^{**} \simeq
{\cp}^{**}$. This implies that the quasibundles $\cp$, ${{\ol
{\cp}}}_{t_0}$, ${{\ol {\cp}}}_{t_1}$ and $\cq$ lie on the same
connected component of the fibre of $h$. This proves the claim.
\begin{flushright} {\it q.e.d} \end{flushright}

\subsubsection{Characteristic classes of principal $H$--bundles}

Let $P$ be a principal $H$--bundle with $H$ a simple algebraic
group. Then by the tables in \cite{bourbaki} for the {\it ``invariant
degrees''}, i.e the degrees of a generating set for the invariant
polynomials under the adjoint action (or equivalently by the Chevalley
restriction theorem, polynomials on the Cartan subalgebra invariant
under the Weyl group action), we see that for $H$--simple, there is a
{\it unique} generator of degree $2$. Since we are over an algebraic
surface, it follows that this invariant polynomial $I_2$ will be the
only one which gives us a characteristic class in $H^{4}(X, \bz)$. We
shall denote this class by $c(P)$. (cf. \cite{beauville})

\begin{exam} {\em When $H$ is classical, then $c(P) = c_2(P(V))$ where $V$ is the
defining representation.}
\end{exam}

We can now state the following key result:

\bprop\label{quasipoly} Let $\cp$ and $\cq$ be two $\mu$--semistable
$H$--quasibundles. Let ${\cp}_0^{**}$ be the canonical {\it
polystable} principal bundle obtained from a choice of associated
graded quasibundle of $\cp$. Let $c({\cp}_0^{**}) = c({\cq}_0^{**})$
for the degree $2$ characteristic class $c$ of the principal bundles
${\cp}_0^{**}$ and ${\cq}_0^{**}$. Then $\cp$ and $\cq$ define the
same point of $M_H(\rho)$ if and only if we have an identification of
pairs $({\cp}_0^{**}, Z({\cp})) \simeq ({\cq}_0^{**}, Z({\cq}))$ where
$Z({\cp})$ is a cycle class in the symmetric power $S^{\it
l_{0}}(X)$ given by
\[
Z({\cp}) := \sum_{x \in X} {length({\ce}_0^{**}/{\ce}_0)_{x} \over
m_{\rho}} \cdot x
\]
with $l_0$ given by:
\[
{\it l}_{0} = {c_2({\ce}_0) - c_2({\ce}_0^{**}) \over m_{\rho}}
\]
and where $\cp$ (resp $\cq$) is obtained from the torsion--free sheaf
$\ce$ (resp $\cf$) by reduction of structure group via $\rho$ and
${\ce}_0$ (resp ${\cf}_0$) is an associated graded sheaf of the
semistable sheaf $\ce$ (resp $\cf$) .\eprop

{\it Proof:} Let us assume that the quasibundles $\cp$ and $\cq$
define the same point in the moduli space $M_H(\rho)$.

We first check that ${\cp}_0^{**} \simeq {\cq}_0^{**}$: suppose that
this does not hold. As we noted in Lemma \ref{polyeq}, it is implies
that for a large degree curve $C$, the restrictions ${\cp}_0^{**}|_{C}
\neq {\cq}_0^{**}|_{C}$. By Lemma \ref{huy}, we can separate the
points corresponding to $\cp$ and $\cq$ in $R_H(\rho)$ by invariant
sections of $\cl$ which contradicts the assumption that they define
the same point in $M_H(\rho)$.

Again, the equality $Z(\cp) = Z(\cq)$ follows immediately since the
torsion--free sheaves $\ce$ and $\cf$ (from which the quasibundles
$\cp$ and $\cq$ are defined) give the same point in the moduli space
$M_G$ by (\cite[Th 8.2.11]{hl}), and the cycle classes $m_{\rho} \cdot
Z(\cp) = Z(\ce)$.

For the converse, suppose that we have an identification:
\[
({\cp}_0^{**}, Z({\cp})) \simeq ({\cq}_0^{**}, Z({\cq})).
\]

By Cor \ref{keycor}, it is enough to check for polystable quasibundles
in the fibre of $h$; by Prop \ref{keycor1} it follows that $\cp$ and
$\cq$ lie in the {\it same connected component of the fibre of} $h:
Q_H({\be},{\it l}) \lr S^{\it l}(X)$.

By the definition of the determinant line bundle on $R_H(\rho)$, it is
obtained by pulling back the corresponding line bundle from the total
family $R_G$. By the theorem of Baranovsky (\cite{bara}) and
Ellingsrud and Lehn (\cite{el}) one knows that the fibres of the
morphism $q: Q({\be},{\it l}) \lr S^{\it l}(X)$ are {\it connected}
and also that the line bundle $\cl$ restricted to this fibre is {\it
  trivial} (cf. \cite[8.2.1]{hl} and Li (\cite[Lemma 3.4,3.5]{li})).

By the commutative diagram seen earlier, we have $q \circ f_{\rho} =
h$. Therefore $\cl$ is trivial on the fibres of $h$ as well. In
particular, the connected components of the fibre of $h$ get mapped to
the same point in $M_H(\rho)$. This proves that the points defined by
$\cp$ and $\cq$ coincide in $M_H(\rho)$. \begin{flushright} {\it
q.e.d}
\end{flushright}

\bprop\label{stableembed} There is a canonical morphism $j : M_H({\sf
  c})^{s} \lr M_H(\rho)$ which is an embedding of the moduli space of
isomorphism classes $\mu$--stable principal $H$--bundles in
$M_H(\rho)$. \eprop

{\it Proof}: The existence of the morphism $j$ follows by the weak
coarse moduli property of the moduli space $M_H(\rho)$. The
injectivity of $j$ follows from Lemma \ref{polyeq}. 

Now choose a general curve $C \in |a\Theta|$ which gives Bogomolov's
effective restriction theorem, it is easily seen that the differential
of the restriction map $r|_{C}$ (by usual arguments involving
Enriques-Severi lemma) is {\it injective}. This immediately implies
that the map $j$ is an {\it embedding}.
\begin{flushright} {\it q.e.d} \end{flushright}

Let ${\goth M}_H$ be defined by the disjoint union:
\[
{\goth M}_H = \coprod_{{\it l} \geq 0} M_H^{{\mu}-poly}({\sf c} - {\it l})
\times S^{{\it l}}(X) \eqno(\ast)
\]
where $M_H^{{\mu}-poly}({\sf c} - {\it l})$ is the moduli space of
$\mu$--semistable principal $H$--bundles with characteristic class
${\sf c} - {\it l}$ (represented as classes of polystable bundles); in
particular, the big stratum is $M_H^{{\mu}-poly}({\sf c}) = M_H^0$.

In conclusion we have the following:

\bth\label{strata} There is a set--theoretic inclusion
\[
{\goth M}_H \subset M_H(\rho)(\bc)
\]
Furthermore, the underlying set of points ${\ol {M_H^{0}}}(\bc)$, of
the closure of the moduli space of equivalence classes of
$\mu$--semistable principal $H$--bundles, is a subset of ${\goth
  M}_H$. In fact, ${\ol {M_H^{0}}}(\bc)$ is independent of the
representation $\rho$ {\sf upto homeomorphism}.  \eeth

{\it Proof:} To complete the proof of the theorem we need only show
that if $({\cp}^{**}, Z)$ is a point on the right hand side of
$(\ast)$ above, then there exists a point in $M_H(\rho)$ which
corresponds to it.

Since ${\cp}^{**}$ is a polystable principal $H$--bundle, using $\rho$
if we extend the structure group to $G = GL(V)$, we get a
$\mu$--polystable locally free sheaf ${\be} \in M_G$ together with a
reduction $\sigma :X \lr {\be}(G/H)$ (as before we confuse the locally
free sheaf $\be$ with the associated principal $GL(V)$--bundle $\be$).

By the results of Li, it is known that the map $j:S^{\it l}(X) \lr
M_G$, induced by the morphism from $Q(\be,{\it l}) \lr M_G$, is an
{\em embedding}. Therefore the element $m_{\rho} \cdot Z \in S^{\it
l}(X)$ gives a polystable torsion--free sheaf $\ce$ such that
${\ce}^{**} \simeq {\be}$ and such that $m_{\rho} \cdot Z = Z(\ce)$.

We therefore have a reduction of structure group of ${\ce}|_{U(\ce)}$
coming from the restriction $\sigma|_{U(\ce)}$ of the reduction of
structure group of $\be$ which gives the principal bundle
${\cp}^{**}$. Viewing this $\sigma|_{U(\ce)}$ as a morphism $U(\ce)
\lr S_{\ce}//H$, we get a canonical extension ${\sigma_1} : X \lr
S_{\ce}//H$. This gives rise to a polystable $H$--quasibundle $\cp \in
M_H(\rho)$ which corresponds to the point $({\cp}^{**}, Z)$. That this
map which sends $({\cp}^{**}, Z) \lr {\cp}$ is injective follows from
the first part of the proof of Prop \ref{quasipoly}.

The last part of the theorem which realises the closure ${\ol
  {M_H^{0}}}(\bc)$ in ${\goth M}_H$ follows from the consequences of
Theorem \ref{dynkimulti} and the definition Def \ref{cycle}. The
uniqueness upto homeomorphism is precisely the content of Prop
\ref{homeo} below.

\begin{flushright} {\it q.e.d} \end{flushright}

\bprop\label{homeo} The closure ${\ol {M_H^{0}}}(\bc)$ is unique upto
homeomorphism. In fact, the isomorphic copies of $M_H^0$ defines a
{\em correspondence} between the closures which gives the
homeomorphism. Furthermore, the normalisations of the closures are
isomorphic as reduced projective schemes.\eprop

{\it Proof}: Since most the points needed in the proof have already
been discussed above, we content ourselves by giving the main steps in
the proof of this Proposition.

Let $\rho_i: H \hra SL(V_i)$ for $i = 1,2$ be two faithful
representations of $H$ and suppose that $P_t(i)$ be family of
$\mu$--semistable principal $H$--bundles with isomorphic generic
fibres and with two limits as {\it polystable $H$--quasibundles}
coming from the representations $\rho_i$. Let the limits be $P(i)$, $i
= 1,2$.

{\it Claim 1}: The double dual polystable $H$--bundles $P(1)^{**}
\simeq P(2)^{**}$. To see this, restrict the families $P_t(i)$ to a
general high degree curve $C$ which avoids the singular loci of $P(i)$
for $i = 1,2$. Then, since for such families over curves the
(polystable) limit is uniquely defined, it follows that
$P(1)^{**}|_{C} \simeq P(2)^{**}|_{C}$. By Lemma \ref{polyeq} the claim
follows.

{\it Claim 2}: The cycles $Z(P(i))$ also coincide. This follows from
two observations, namely that for the limits the open subsets where
the $P(i)$ are genuine principal bundles coincide and hence the sets
$Sing(P(i))_{red}$ coincide. The multiplicities at all the singular
points also coincide by the discussion before Definition \ref{cycle}.

These two claims imply the set--theoretic identification of the
closures. Moreover, the above discussion gives an identification of
closures of any curve in $M_H^0$. This implies that the projections to
$M_H(\rho_i)$ from the graph $\Gamma \subset M_H(\rho_1) \times
M_H(\rho_2)$, which is closed in the product, actually give
homeomorphisms.  The comment on normalisations now follows by
Zariski's main theorem.

\bcor\label{donu} The closure of subset $M_H({\sf c})^{s} \subset
M_H^0$ gives the Donaldson-Uhlenbeck compactification of the moduli
space of principal $H$--bundles with irreducible ASD connections. 
\ecor

\subsubsection{Relationship with G\'omez-Sols and Schmitt's moduli space}

We work in the set-up of Schmitt's recent paper (cf. \cite{schmitt} and
\cite{schmitt1}). Fix a faithful representation $\rho: H \hra
SL(V)$. Define the moduli functors
\[
{\underline M}(\rho)^{ss}_{P}: {\underline {Sch}}_{\bc} \lr {\underline{Set}}
\]
which sends $S$ to ``equivalence classes of families of
Gieseker--Maruyama semistable honest singular principal $H$--bundles
(or $H$--quasibundles) with Hilbert polynomial P on $X$ parametrised
by $S$''. Then the main theorem of \cite{schmitt1} is that there is a
projective scheme ${\cm}(\rho)^{ss}_P$ which {\it coarsely represents}
the functor ${\underline M}(\rho)^{ss}_{P}$. The notions of
Gieseker--Maruyama semistability of quasibundles is defined in
\cite{gomez} and \cite{schmitt1} and it is also shown that if ${\cp} =
(\ce,\tau)$ is an $H$--quasibundle then the following equivalence
holds (\cite[Section 5.1]{schmitt1}): ${\cp}~~is~~{\mu}-stable
\implies$ ${\cp}~~is~~Gieseker-Maruyama-stable \implies$
${\cp}~~is~~Gieseker-Maruyama-semistable \implies$
${\cp}~~is~~{\mu}-semistable$. These implications, together with the
coarse moduli property of the functor ${\underline M}(\rho)^{ss}_{P}$
implies that there is a morphism:
\[
{\cm}(\rho)^{ss}_P \lr M_H(\rho).
\]

At the risk of repetition, we now summarise the above results in the
following:

\bth\label{uhlenbeck} Let $H$ be a semisimple algebraic group. 
\begin{enumerate}
\item There exists a projective scheme $M_H(\rho)$ which parametrises
equivalence classes of $\mu$--semistable $H$--quasibundles with
fixed characteristic classes, on the smooth projective surface
$X$. 

\item This has an open subscheme of equivalence classes
$\mu$--semistable principal $H$--bundles $M_H^{0}$.

\item If $M_H^s$ is the subscheme of $M_H(\rho)$ consisting of the
stable principal $H$--bundles, then its closure $\overline {M^s}$ in
$M_H(\rho)$ corresponds to the Donaldson-Uhlenbeck compactification of
the moduli space of $\mu$--stable vector bundles on $X$ (or
equivalently of anti self--dual connections).

\item {\it When $H$ is a simple algebraic group} the points of the
moduli space ${\ol {M_H^{0}}}$ are given by the pairs
$({\cp}_0^{**},Z(\cp))$ where, to each semistable $H$--quasi bundle
$\cp$, we associate ${\cp}_0^{**}$ the canonical polystable principal
$H$--bundle obtained from the quasibundle $\cp$ and the cycle class
$Z(\cp)$ which lies in $S^{l_0}(X)$, ${l_0}$ being given by:
\[
{c_2({\ce}_0) - c_2({\ce}_0^{**}) \over m_{\rho}}.
\]
Note that $\cp$ is obtained from the torsion--free sheaf $\ce$ by
reduction of structure group via $\rho$ and ${\ce}_0 = gr^{\mu}(\ce)$
and $m_{\rho}$ is the Dynkin index associated to $\rho$.

\item {\it When $H$ is a simple algebraic group} the underlying set of
  points of the moduli space ${\ol {M_H^{0}}}$, upto homeomorphism, is
  independent of the representation $\rho$.

\item By (\cite[Main theorem]{schmitt1}), we have a morphism
$${\cm}(\rho)^{ss}_P \lr M_H(\rho)$$. 

\item In a certain sense, ${\ol {M_H^{0}}}$ is the minimal
compactification of $M_H^{0}$.

\end{enumerate}

\eeth

\section{Non--emptiness of the moduli space}

In this section $H$ is an arbitrary semisimple algebraic group. The
aim of this section is to prove that the moduli space $M_{H}(\sf
c)^{s}$ of $\mu$--stable principal $H$--bundles on a smooth projective
surface $X$ for {\it large} characteristic classes $\sf c$ is
non--empty. When $H = SL(2)$ this is a theorem due to Taubes
\cite{taubes} and later due to Gieseker \cite{gies}. We quote their
result:

\bth\label{giestaub} (Taubes, Gieseker) For any choice of polarisation
$\Theta$, for all $c \geq 2p_g(X) + 2$, the moduli space
$M_{_{SL(2)}}(c)^s$ of stable rank $2$ bundles with trivial
determinant and $c_2 = c$ is non--empty. \eeth

\subsection{\it Stable bundles and the principal three dimensional subgroup}

We begin by recalling a few facts about semisimple Lie algebras. Let
$\goth {H} = Lie(H)$. Then one knows that there exist the so--called
{\it principal} ${\goth {sl(2)}}$'s inside ${\goth H}$ which are
distinguished by the property that they do not lie in any parabolic
subalgebra of $\goth H$. This gives a representation $\phi: SL(2)
\rightarrow H$ which is {\it irreducible}, in the sense that
$Im(\phi)$ is not contained in any proper parabolic subgroup of $H$.
(cf. \cite{bourbaki} Chapter 8, Exercise 5 \S11, page 246).

We fix one such homomorphism of a principal $SL(2) \rightarrow H$.

Let $C$ be a smooth projective curve of genus $g \geq 2$. Let $V$ be a
stable vector bundle of rank $2$. Then by the Narasimhan--Seshadri
theorem one knows that there exists an irreducible representation
$\rho: {\pi}_1(C) \lr SU(2)$ such that $V \simeq V_{\rho}$.

\bdefe We define the monodromy subgroup ${\mathcal M}(V)$ of $V \simeq
V_{\rho}$ to be the Zariski closure of $Im(\rho)$ in $SL(2)$. It can
be viewed as the minimal subgroup to which the structure group of $V$
can be reduced.  \edefe

\brem Recall that one knows that the {\it reductivity} of ${\mathcal M}(V)$ is
equivalent to the polystability of $V$.  Since one knows  
the set of all finite subgroups of $SL(2)$ (this is classical, for example cf.
\cite{serre1}), the only
possibilities for ${\mathcal M}(V)$ are the following:
\begin{enumerate}
\item Finite cyclic groups $C_n$ and the dihedral groups $D_n$.

\item The alternating groups $A_4$ and $A_5$ and the permutation group $S_5$.

\item The whole of $SL(2)$ or its maximal torus.

\end{enumerate}

Of these, since $V$ is stable, we can omit the maximal torus and the
cyclic groups $C_n$. So ${\mathcal M}(V)$ can either be the
alternating groups, $S_5$, or the dihedral groups apart from the whole
of $SL(2)$.

We wish to estimate the set $Z_C$ of representations of ${\pi}_1(C)$ in
$SL(2)$ which lie entirely in these families of finite groups upto
conjugacy by the diagonal action of $SL(2)$.

It is not hard to see that since ${\pi}_1(C)$ is given by $2g$
generators with a single relation, this set $Z$ is at most {\it
countable}. This implies the following lemma.

\erem
 
\blem The locus of points in the moduli space of stable vector bundles
$M_{C,SL(2)}$ of rank $2$ and trivial determinant on the curve $C$
whose monodromy is among the set of finite subgroups listed above has
countable cardinality. \elem

Using this we have the following:

\bprop\label{principal} There exists a rank $2$ stable bundle $E$ with
$c_2(E) \gg 0$ and trivial determinant on the surface $X$ such that
the restriction $E|_{C}$ to general curve $C \subset X$ of high degree
has monodromy subgroup to be the whole of $SL(2)$ itself. \eprop

{\it Proof:} For $c = c_2(E) \gg 0$ (by Theorem \ref{giestaub}), it is
known that the moduli space $M_X(2, {\co}, c_2)^{s} = M_{_{SL(2)}}(c)^s$
is non--empty and has dimension $d(c) = 4 c - 3 {\chi}(\co_X)$. Since
$c_2(E) \gg 0$, it follows that we can make $d \geq 1$.

Now we consider the restriction map $r_C : M_{_{SL(2)}}(c)^{s} \lr
M_{_{C,SL(2)}}^{s}$. If the curve $C$ is chosen sufficiently general
and high degree $k$ then by the Mehta--Ramanathan theorem and
\cite[Prop 2.2]{friemorg}, it can be seen that for $k \gg 0$, the map
$r_C$ is in fact an {\it immersion}. This implies that $Im(r_C)
\subset M_{_{C,SL(2)}}^{s}$ is not completely contained in the subset
$Z_C$ of stable bundles with finite monodromy groups. Therefore there
exists at least one stable bundle with the whole of $SL(2)$ as its
monodromy subgroup.  This proves the proposition.

Before we proceed to our next lemma we recall the following
formulation of {\it irreducible} sets (\cite[page 129]{rama}):

\bdefe Let $K$ be a maximal compact subgroup of $H$. A subset $A
\subset K$ is said to be {\em irreducible} if
\[
\{ Y \in {\goth H} | ad~x(Y) = Y, ~~\forall~~ x \in A \} =
centre~of~{\goth H}.
\]
A (homomorphism) representation $\delta: \Gamma \lr K$ is said to be
irreducible if the image ${\delta}(\Gamma)$ is irreducible.\edefe

Note that in our case since $H$ is semisimple $centre~of~{\goth H}$ is
trivial. We have the following equivalences and we quote:
\bprop\label{ramanath}(\cite[Prop 2.1]{rama}) Let $K$ be a maximal
compact subgroup of $H$ and let $A \subset K$ be a subgroup. Then the
following are equivalent
\begin{enumerate}
\item $A$ is irreducible.
\item $A$ leaves no parabolic subalgebra of ${\goth H}$ invariant.
\item For any parabolic subgroup $P \subset H$, $A$ acts without fixed
points on $H/P$.
\end{enumerate}
\eprop

\blem\label{subbu} If $V$ is a stable bundle on a curve $C$ such that
${\mathcal M}(V) = SL(2)$ then $V(\phi)$ is a stable $H$--bundle on
$C$ by an extension of structure group via the principal homomorphism
$\phi$. \elem

{\it Proof:} Let $\rho: {\pi}_C \lr SU(2)$ is an irreducible unitary
representation and let $V \simeq V_{\rho}$ be a stable bundle on $C$
such that ${\mathcal M}(V) = SL(2)$.

This implies that $Im(\rho)$ is Zariski dense in $SL(2)$. Since
$Im(\phi)$ is not completely contained in any parabolic subgroup of
$H$ by Prop \ref{ramanath} we have the following:
\[
\{ Y \in {\goth H} | ad~x(Y) = Y, ~~\forall~~ x \in {\goth {sl}(2)} \} =
centre~of~{\goth H} = trivial.
\]
In other words, we can say that $Im(\phi)$ is an {\it irreducible
subset} in $H$ (or equivalently, by the ``unitarian trick'',
$Im(\phi|_{\goth {su}(2)})$ is an {\it irreducible} subset of $K$,
where as above, $K \subset H$ is a maximal compact subgroup).

By \cite[Prop 2.2]{rama}, to show that $V(\phi)$ is stable as an
$H$--bundle, we need only show that the representation 
\[
\eta = (\phi \circ \rho): {\pi}_1(C) \lr K
\]
is {\it irreducible}.

If $ Y \in {\goth H}$ is such that $ad~x(Y) = Y, ~~\forall~~ x \in
Lie(Im (\eta)) $, then by the density of $Im(\eta)$ in $Im(\phi) =
{\goth {sl}(2)}$ and hence by continuity, we see that for such an
element $Y$, one has $ad x(Y) = Y, ~~\forall~~ x \in {\goth {sl}(2)}$
and hence $Y \in centre~of~{\goth H}$. This shows that $Im(\eta)$ is
an irreducible set in $K \subset H$ and we are done. ~~~{\it q.e.d}

We can now conclude the following:

\bth\label{taubes} The moduli space $M_H(c)^{s}$ of $\mu$--stable
principal $H$--bundles on a smooth projective surface $X$ is
non--empty for $c > \delta$, where $delta$ depends only of $p_g(X)$
and not on the polarisation $\Theta$ on $X$. \eeth

{\it Proof:} We {\it claim} that if $\phi : SL(2) \lr H$ is as above a
principal $SL(2)$ in $H$ then any $E$ as in Prop \ref{principal} has
the property that $E(\phi)$ is a {\it stable} principal $H$--bundle.

By the converse to the Mehta--Ramanathan theorem (Lemma \ref{rest}),
we see that it is enough to prove that $E(\phi)|_{C}$ is {\it
stable}. Since ${\mathcal M}(E|_{C}) = SL(2)$ this is immediate by
Lemma \ref{subbu}. The {\it largeness} of the characteristic classes
of the associated $H$--bundle is determined by the largeness of the
$c_2$ of the rank $2$ bundle $E$. This can determined by the general
methods of Borel and Hirzebruch (see the recent preprint Beauville
\cite{beauville} for this).
\begin{flushright} {\it q.e.d} \end{flushright}

\subsection{Concluding remarks}
 
\begin{enumerate}
\item The questions addressed in this paper can be posed in positive
characteristic as well, but as in the case of curves a subtler
analysis of representation theoretic bounds such as the ones
considered in \cite{bapa} will have to be carried out. We hope to do
this in future. 

\item The basic questions regarding irreducibility,
reducedness, generic smoothness and normality are yet to be answered
for these moduli spaces.

\end{enumerate}

%\nopagebreak[3]


\begin{thebibliography}{99999}
  
\bibitem{atiyahbott} M.F. Atiyah and R. Bott: The Yang-Mills equations
  over Riemann Surfaces, {\it Phil. Tr. R. Soc. Lond. A} {\bf 308}
  (1982), 523-615.


\bibitem{atiyah} M.Atiyah, N.Hitchin and I.M.Singer : Self--duality in
four dimensions, {\it Proc. Royal Soc. Lond} A {\bf 362} (1978), 425-461. 

\bibitem{basa} V.Balaji and C.S.Seshadri: Semistable Principal
Bundles-I (in characteristic zero), {\it Journal of Algebra} {\bf 258}
(2002), 321-347.

  
\bibitem{bapa} V.Balaji and A.J.Parameswaran : Semistable principal
  bundles-II (in positive characteristics).  {\it Transformation
  Groups}, {\bf l8}, No 1, (2003), pp 3-36.


\bibitem{bara} V.Baranovsky : Moduli of sheaves on surfaces and action
  of the oscillator algebra, {\it J. Diff. Geom}, {\bf 55}, No 2,
  (2000), pp 193-227.


\bibitem{beauville} A.Beauville : Chern classes for representations of
reductive groups, math.AG/0104031.


\bibitem{bisgom} I.Biswas and T.L.G\'omez : Restriction theorem for
principal bundles, {\it Math.Annalen}, {\bf 327}, (2003), 773-792. 

\bibitem{bourbaki} N.Bourbaki : Groupes et alg\`ebres de Lie, Chap
7-8, Paris, Hermann, (1975).

\bibitem{brave} A.Braverman, M.Finkelberg and D.Gaitsgory : Uhlenbeck spaces 
via affine Lie algebras, math.AG/0301176.


\bibitem{collio} J-L Colliot-Th\'el\`ene and J.Sansuc : Fibr\'es
  quadratiques et composantes connexes r\'eelles, {\it Math.Annalen},
  {\bf 244} (1979) pp 105-134.


\bibitem{donald} S.K.Donaldson : Polynomial invariants for smooth four
manifolds, {\it Topology} {\bf 29}, (1990), 257-315.

\bibitem{dk} S.K.Donaldson and P.B.Kronheimer: The Geometry of
  Four-Manifolds, Oxford (1990).

\bibitem{dynkin} E.B.Dynkin : Semisimple subalgebras of semisimple Lie
algebras, {\it Am Math Soc Transl} Ser II {\bf 6},(1957) 111-244.


\bibitem{el} G.Ellingsrud and M.Lehn : Irreducibility of the punctual
  quotient scheme of a surface. {\it Ark. Mat. 37} (1999), no. 2,
  245--254.

%\bibitem{faltings} G.Faltings

\bibitem{friemorg} R.Friedman and J.Morgan : Smooth four--manifolds
and complex surfaces, Ergebnisse, Vol 27, Springer Verlag (1994).

\bibitem{gies} D.Gieseker : A construction of stable bundles on an
algebraic surface, {\it J. Diff. Geom} {\bf 27}, (1988), 137-154.

\bibitem{giesli} D.Gieseker and J.Li : Moduli of high rank 
vector bundles over surfaces, {\it J.AMS} {\bf 9}, (1996), 107-151. 

\bibitem{gomez} T.G\'omez and I.Sols : Moduli spaces of principal
sheaves over projective varieties, {\it Annals of Math.} (to appear).

%\bibitem{dmil} P.Deligne and J.Milne: Tannaka Categories, Springer
%Lecture Notes in Mathematics Vol 900.

\bibitem{hl} D.Huybrechts and M.Lehn: The Geometry of moduli spaces of
sheaves, Aspects of Mathematics E31, Vieweg, Braunschweig/Wiesbaden,
(1997).

\bibitem{langer} A.Langer: Semistable sheaves in positive
  characteristic, {\it Annals of Math.}

\bibitem{snr} S.Kumar, M.S.Narasimhan and A.Ramanathan, Infinite
Grassmanians and moduli spaces of $G$--bundles, {\it Math Ann} {\bf
300}(1994), 41-75.

\bibitem{langton} S.Langton: Valuative criterion for families of vector
bundles on algebraic varieties, {\it Annals of Mathematics}(2) {\bf
101} (1975) pp 88-110.

\bibitem{lepot} J.Le Potier: Fibr\'e d\'eterminant et courbes de saut
sur les surfaces alg\'ebriques, {\it Complex projective Geometry},
London Mathematical Society : Bergen(1989), 213-240.

\bibitem{li} J.Li: Algebraic geometric interpretation of Donaldson's
polynomial invariants of algebraic surfaces, {\it J.Diff.Geom} {\bf
37} (1993), 416-466.

\bibitem{li1} J.Li:Compactification of moduli of vector bundles over
  algebraic surfaces, (In a volume of papers dedicated to Prof C.G.Gu).
  
\bibitem{maruyama} M.Maruyama: The Theorem of
  Grauert-M\"ulich-Spindler, {\it Math. Ann } {\bf 255},(1981)
  317--333.

%\bibitem{maruyama2} M.Maruyama: Moduli of stable sheaves II,
 % {\it J.Math.Kyoto Univ.} {\bf 18} (1978), 557-614.

\bibitem{mr} V.B.Mehta and A.Ramanathan: Restriction of stable sheaves
and representations of the fundamental group, {\it Invent.Math} {\bf
77}, (1984), 163--172.

\bibitem{morgan} J.Morgan : Comparison of the Donaldson polynomial
invariants with their algebro geometric analogues, {\it Topology} {\bf
32} (1993), 449-488.

\bibitem{ns} M.S. Narasimhan and C.S.Seshadri: Stable and unitary
vector bundles on a compact Riemann surface, {\it Annals of
Mathematics }(2){\bf 82} (1965) pp 540-567.


\bibitem{newstead} P.Newstead: Introduction to moduli problems, TIFR
  Lecture Notes in Mathematics, (1978).

%\bibitem{nori} M. V. Nori: The fundamental group scheme,
%{\it Proc.Ind.Acad.Sci (Math.Sci)} {\bf 91} (1982), 73--122.



%\bibitem{okonek} Okonek, Schneider and Spindler : Vector bundles on
%projective spaces. {\it Birkhauser Verlag}

\bibitem{rama} A.Ramanathan: Stable principal bundles on a compact
  Riemann surface, {\it Math.Ann}, {\bf 213}, (1975) 129-152.

\bibitem{r1} A. Ramanathan: Stable principal bundles on a compact
  Riemann surface - Construction of moduli space ({\it Thesis}, Bombay
  University 1976) {\it Proc.Ind.Acad.Sci} {\bf 106}, (1996), 301-328
  and 421-449.
  
  
\bibitem{rs} A.Ramanathan and S.Subramaniam : Einstein--Hermitian
connections on principal bundles and stability, {\it J.Reine
Angew. Math} {\bf 390} (1988), 21-31.


\bibitem{rr} S.Ramanan and A.Ramanathan : Some remarks on the
instability flag, {\it Tohoku Math. Journ} {\bf 36} (1984), 269-291.


\bibitem{schmitt} A.Schmitt: Singular Principal Bundles over higher
dimensional manifolds and their moduli spaces, {\it IMRN}, {\bf 23}
(2002), 1183-1209.

\bibitem{schmitt1} A.Schmitt: A closer look at semistability for
singular principal bundles, (to appear in {\it IMRN}).


\bibitem{serre} J.P.Serre: Espaces fibr\'es alg\'ebriques, Anneaux de
Chow et applications, S\'eminaire Chevalley, (1958), also {\it
Documents Math\'ematiques} {\bf 1} (2001), 107-140.

 
\bibitem{serre1} J.P.Serre: Moursund Lectures, University of Oregon
Mathematics Department (1998).

\bibitem{sim} C. Simpson: Higgs bundles and Local systems,
{\it Pub. I.H.E.S.} {\bf 75} (1992), 5-95.

\bibitem{taubes} C.H.Taubes : Self-dual connections on $4$--manifolds
with indefinite intersection matrix, {\it J.Diff.Geom} {\bf 19}
(1984), 517-560.









\end{thebibliography}
\end{document}